\documentclass[a4paper,10pt]{article}
\usepackage[utf8]{inputenc}
\clearpage{}\usepackage{graphics}
\usepackage{amsmath}
\usepackage{amssymb}
\usepackage{amsthm}
\usepackage{enumerate}
\usepackage{siunitx}
\usepackage{stmaryrd} \usepackage{dsfont} \usepackage{graphicx}
\usepackage{enumerate}
\usepackage{subfigure}
\usepackage{doi}
\usepackage{cite}

\usepackage{mdframed}

\usepackage{acronym}
\usepackage{mhchem}

\usepackage{xcolor}

\newcommand{\change}[1]{#1} 

\newtheorem{lemma}{Lemma}[]

\newcommand{\ddt}[1]{\frac{d #1}{dt}}

\newcommand{\DDT}[1]{\frac{D #1}{Dt}}
\newcommand{\inproduct}[2]{\left\langle #1 , #2 \right\rangle}

\DeclareMathOperator{\sgn}{sgn}
\DeclareMathOperator{\gr}{gr} \newcommand{\transpose}{{\sf T}}

\newcommand{\RR}{\mathds{R}}  

\newcommand{\nsigma}{n_\Sigma}

\newcommand{\ndomega}{n_{\partial\Omega}}

  \newcommand{\normalspeed}{V_\Sigma} 

\newcommand{\phistar}{\phi}
\clearpage{}

\usepackage{amsthm}
\usepackage[a4paper,left=2.5cm,right=4cm,top=3.5cm,bottom=5cm]{geometry}
\usepackage{authblk} 

\usepackage[margin=1cm]{caption} \usepackage{cite}

\date{}

\theoremstyle{definition}

\newtheorem{remark}{Remark}

\title{A locally signed-distance preserving level set method (SDPLS) for moving interfaces}

\author[1]{Mathis Fricke\thanks{fricke@mma.tu-darmstadt.de (Corresponding Author)}}
\author[1]{Tomislav Marić\thanks{maric@mma.tu-darmstadt.de}}
\author[1]{Aleksandar Vučković \thanks{aleksandar.vuckovic@stud.tu-darmstadt.de}}
\author[2]{Ilia Roisman \thanks{roisman@sla.tu-darmstadt.de}}
\author[1]{Dieter Bothe\thanks{bothe@mma.tu-darmstadt.de}}
\affil[1]{Mathematical Modeling and Analysis, TU Darmstadt}
\affil[2]{Institute for Fluid Mechanics and Aerodynamics, TU Darmstadt}

\begin{document}

\maketitle

\abstract{
It is well-known that the standard level set advection equation does not preserve the signed distance property, which is a desirable property for the level set function representing a moving interface. Therefore, reinitialization or redistancing methods are frequently applied to restore the signed distance property while keeping the zero-contour fixed. As an alternative approach to these methods, we introduce a modified level set advection equation that intrinsically preserves the norm of the gradient at the interface, i.e.\ the local signed distance property. Mathematically, this is achieved by introducing a carefully chosen source term being proportional to the local rate of interfacial area generation. The introduction of the source term turns the problem into a non-linear one. However, we show that by discretizing the source term explicitly in time, it is sufficient to solve a linear equation in each time step. Notably, without further adjustment, the method works in the case of a moving contact line. This is a major advantage since redistancing is known to be an issue when contact lines are involved. We provide a first implementation of the method in a simple first-order upwind scheme in both two and three spatial dimensions.
} \newline

\section{Introduction}
Since the beginning of computational fluid dynamics in the 1960s, researchers in the field of computational multiphase flows have developed a variety of powerful specialized methods that allow the numerical simulation of rather complex multiphase flow problems. For example, modern multiphase flow methods allow to simulate flows with strong deformations and topological changes such as breakup and coalescence of droplets and the atomization of liquid jets. Moreover, additional transport processes are studied such as evaporation and condensation or the transport of surface-active substances leading to local Marangoni forces. Even reactive mass transfer processes involving dissolved chemical species are accessible with numerical simulation. Clearly, there is no single method that is equally suitable for all possible flow configurations. Instead, there is a variety of numerical methods with specific advantages and disadvantages for a given multiphase flow problem. Here, we only mention these different methods briefly. For a comprehensive review of numerical methods in multiphase flows, the reader is referred to \cite{Scardovelli1999,Tryggvason2011,Gross2011} and the references given therein.\\
The Volume-of-Fluid (VOF) method introduced by Hirt and Nichols \cite{Hirt1981} is designed to be volume conservative even in the discrete case. This is achieved by transporting the phase indicator function with a finite volume discretization. However, in practice, exact volume conservation may not be achieved (e.g., because of directional splitting and heuristic volume redistribution algorithms to ensure boundedness). The VOF method requires some specialized techniques to numerically transport the sharp discontinuity of the phase indicator function. Geometrical VOF methods (see, e.g., \cite{Tryggvason2011,Maric2020} for an overview) reconstruct the sharp interface after each transport step in order to prevent numerical diffusion of the interface. In general, the numerical approximation of the interface curvature is known to be challenging for VOF methods, in particular for unstructured meshes \cite{Popinet2018,Maric2020}. 
Another important class of methods are front tracking methods (see, e.g, \cite{Unverdi1992}). They represent the interface by a set of connected marker points which are advected by the flow in a Lagrangian manner. While the accuracy of these methods is potentially very high and the interface stays sharp by construction, there is no built-in phase volume conservation and handling topological changes is expensive.\\
\\
In the present work, we focus on the level set method due to Osher and Sethian \cite{Osher1988,Sethian2001,Osher2003,Sethian2003,Gibou2018,Giga.2006}. It represents the interface $\Sigma(t)$ as the zero contour of a smooth function $\phi$ called the ``level set function'' according to
\[ \Sigma(t) = \{ x \in \Omega: \phi(t,x) = 0 \}. \]
Similarly to the VOF method, the level set method is able to handle topological changes naturally. The higher regularity of the level set function compared to the phase indicator function in VOF has some numerical advantages. For example, the interface curvature can be computed more accurately at a given numerical effort. On the other hand, the level set method is not strictly volume conservative in the discrete case. This is a consequence of the fact that the level set $\phi$ is an auxiliary function that has no physical meaning except for the zero iso-surface representing the interface. Even though the integral of $\phi$ over the computational domain is conserved\footnote{Here, for simplicity, we assume that there is no in- or outflow at the domain boundary $\partial\Omega$.\newline Then $\partial_t \phi + v \cdot \nabla \phi = 0$ implies $\ddt{} \int_\Omega \phi \, dV=0$.} for an incompressible flow, one cannot obtain relevant physical information from this conservation principle in general. If, however, a smoothed Heaviside function is used, this conservation principle becomes the basis for the conservative level set method due to Olsson and Kreiss \cite{Olsson2005,Olsson2007,Zahedi2009} (see below). \\
\\
A prominent example of a level set function is the so-called signed distance function $d_\Sigma$ associated with a given hypersurface $\Sigma$.  It is characterized by the property $|\nabla d_\Sigma| = 1$ and has a number of convenient properties (see Appendix~\ref{sec:appendix_interfaces_and_kinematics} for more details). For example, the mean curvature and normal speed of $\Sigma$ be computed via the Laplacian and the partial time-derivative, respectively, i.e.\
\[ \kappa_\Sigma = \Delta d_\Sigma, \quad \normalspeed = - \partial_t d_\Sigma \quad \text{on} \ \Sigma. \]
By keeping (at least approximately) the signed distance property, one can assure that the gradients of the level set function become neither too steep nor too flat compared to the computational resolution. This is important for numerical accuracy and stability of the method. But, as is well-known, the classical level set equation
\begin{align}\label{eqn:level_set_equation}
\partial_t \phi + v \cdot \nabla \phi = 0
\end{align}
does not preserve the signed distance property, i.e.\ the norm of the gradient of $\phi$ will not be constant, not even for the exact solution of \eqref{eqn:level_set_equation} and not even at $\phi=0$. This well-known problem is addressed by ``reinitialization'' or ``redistancing'' methods that compute a signed distance field from a given fixed zero contour. We refer to the references \cite{Sussman1994,Sethian1996,Sussman1999} for some of the first contributions in this direction. The recomputed level set is then used to replace the ``degenerated'' level set function. In many cases, the redistancing is done by solving a PDE in pseudo time. For example, Sussman et al.\ \cite{Sussman1994} proposed to solve the initial value problem for the Hamilton-Jacobi equation
\begin{align}\label{eqn:reinit_sussman1994}
\frac{\partial}{\partial \tau} \, \phi + \sgn(\phi^0) (|\nabla \phi| - 1) = 0, \quad \phi|_{\tau=0} = \phi^0
\end{align}
to transform a given level set function $\phi^0$ into a signed distance field while keeping the zero contour fixed. In order to save computational time, the redistancing is usually only applied on demand and not at every time step.\\
\\
There is a particular issue with the redistancing algorithms when the free surface touches the boundary to form a so-called ``contact line''. Physically, the contact line (defined as the line of intersection of the free surface with a solid boundary) plays an important role because its presence leads to a variety of additional physical effects that are to be modeled via appropriate boundary conditions. Leaving aside the details of the mathematical modeling, the redistancing algorithm based on \eqref{eqn:reinit_sussman1994} faces the problem that additional boundary conditions are required in a region close to the contact line called the ``blind spot'' \cite{Rocca2014,Pertant2021}. There is no obvious choice for these boundary conditions because it would require a continuation of the zero level set beyond the domain boundary. Della Rocca and Blanquart reported in \cite[p.35]{Rocca2014} that there are several methods to populate the blind spot which in many cases cause poor calculation of interface curvature and thereby lead to unintended parasitic currents. Moreover, the redistancing itself may contribute to errors in phase volume conservation by an unintended motion of the interface \cite{Osher2003,Solomenko2017}.\\
Olsson and Kreiss addressed the problem of phase volume conservation with a level set method \cite{Olsson2005,Olsson2007,Zahedi2009} that aims at transporting a smoothed step function instead of a signed distance function. A compression scheme in pseudo-time is applied in an intermediate step to restore the characteristic width of the transition region after advection. They showed that this strategy is able to improve the phase volume conservation significantly.\\
\\
There are already some methods in the literature which aim to combine the advection and redistancing of the level set into one single ``monolithic'' approach \cite{Li2005,Li2006,Toure2016,Ouazzi2018,Ville2011,Coupez2007,Bahbah2019,Grave2020,Bonito2015,Luna2019}. Li et al.\ \cite{Li2005,Li2006} and later Touré and Soulaïmani \cite{Toure2016} \change{and Ouazzi et al.\ \cite{Ouazzi2018} applied} a variational formulation where deviations from a signed distance function are penalized in the objective functional. The evolution equation of the level set is then obtained as the corresponding gradient flow to minimize the objective functional. This way the level set function is kept close to a signed distance function during its evolution. However, a displacement of the zero contour with respect to standard level set equation cannot be ruled out completely but is numerically reduced using another penalty term \cite{Luna2019}. Coupez et al.\ \cite{Coupez2007,Ville2011,Bahbah2019,Grave2020,Bonito2015} introduced the ``convected level set method'' that combines the redistancing equation \eqref{eqn:reinit_sussman1994} with the usual advection equation into one single equation. The resulting PDE aims at keeping the level set function close to a signed distance while the motion of the zero contour remains exactly the same on the continuous level (up to discretization errors). Recently, de Luna et al.\ \cite{Luna2019} introduced a conservative monolithic level set method. Similarly to the approach by Li et al.\ \cite{Li2005,Li2006}, an approximate signed distance function is kept by penalizing deviations from a signed distance. The conservation property is achieved by transporting a smoothed step function similarly to the method by Olsson and Kreiss. Finally, it is also worth mentioning that there have been various developments on hybrid methods that combine ideas and methods from different approaches. The hybrid Level Set / Front Tracking method \cite{Shin2011} resolves the challenge of topological interface changes efficiently by tracking the fluid interface as a zero level-set from a signed distance, recovered from the geometrical approximation of the evolving interface in its near vicinity. This way, the interface can be efficiently recovered as a zero-level set in the case of merging or coalescence.\\
\\
In this work, we introduce a modified transport equation for the level set function that ensures exact conservation of $|\nabla \phi|$ at the zero contour on the continuous level. This is done by introducing a non-linear source term that is derived from the kinematics of interfaces as studied in \cite{Fricke2019,Fricke2021,Fricke2018}. We show that, in practice, it is sufficient to solve a linear problem in each timestep. Moreover, the method works for a moving contact line without special treatment. In particular, we demonstrate numerically that the transported contact angle and the curvature at the contact line are converging with mesh refinement.\\
\\
The remainder of this article is organized as follows: The derivation of the level set equation with source term is presented in Section~\ref{sec:eqn_derivation}. The resulting PDE is then discretized with a simple first-order upwind method in Section~\ref{sec:upwind}. First numerical results are discussed in Section~\ref{sec:results}.

\section{A modified level set advection equation}\label{sec:eqn_derivation}
\subsection{Problem setup} 
\begin{figure}[ht]
\centering
\subfigure[Setup without contact line.]{\includegraphics[width=0.34\columnwidth]{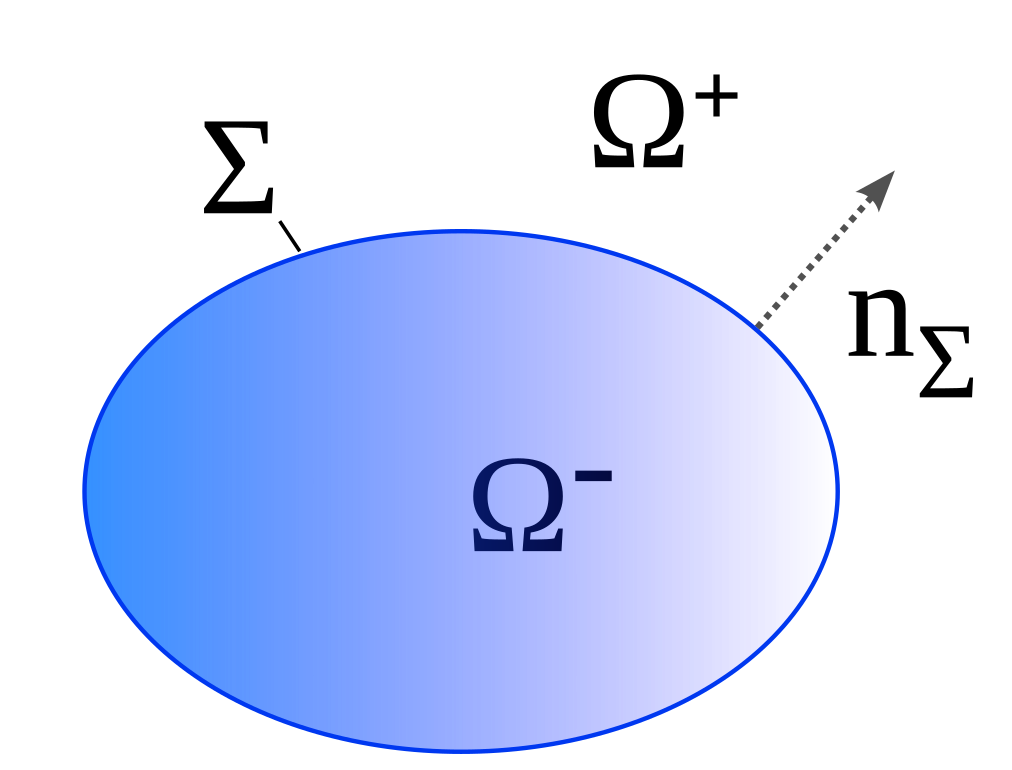}}
\subfigure[Setup with contact line.]{\includegraphics[width=0.55\columnwidth]{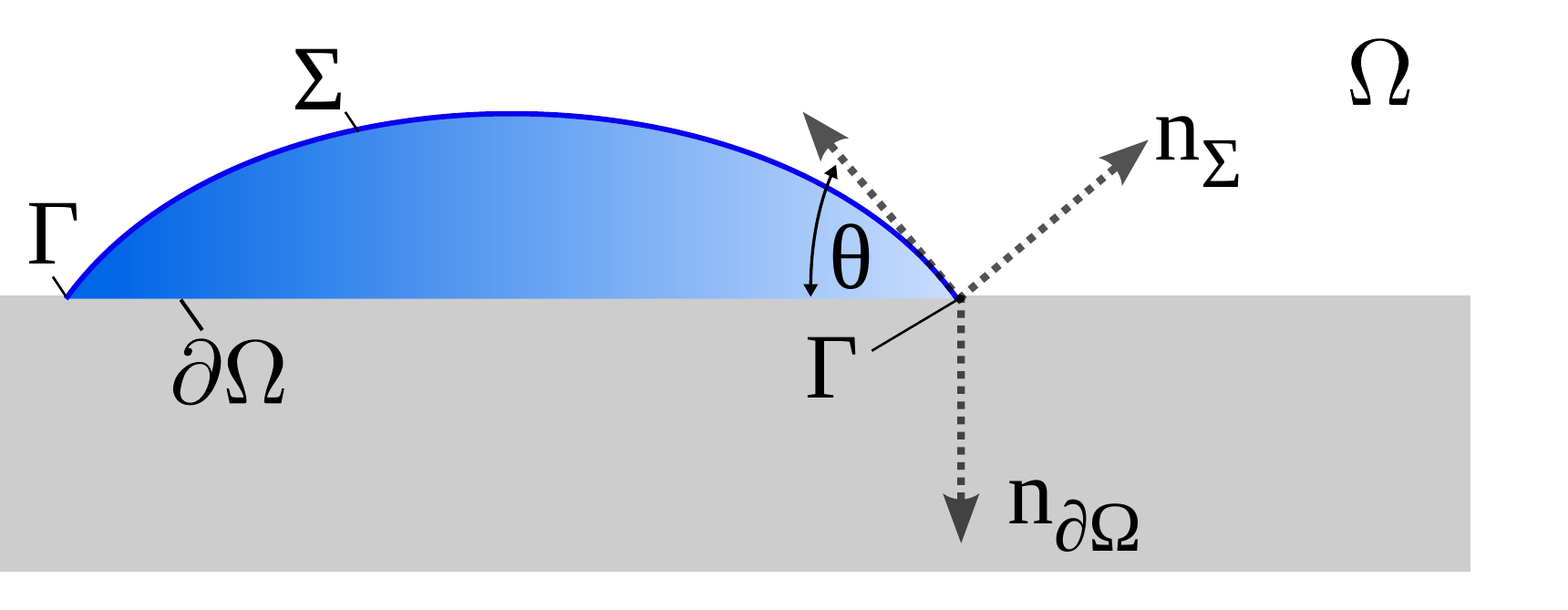}}
\caption{Notation and problem setup.}\label{fig:notation}
\end{figure}
In the present article, we focus on the problem of \emph{advecting} the fluid interface in a two-phase flow problem. We assume the transporting velocity field as given while keeping in mind that, in applications of the method, it will emerge as a (numerical) solution of the Navier Stokes equations. Mathematically, we consider a domain $\Omega$ with piecewise smooth boundary 
$\partial\Omega$ and velocity field $v=v(t,x)$ that is of class $\mathcal{C}^1$ on $\RR \times \Omega$ and divergence-free (modeling incompressible flow)\footnote{Note that in a real two-phase flow problem, the velocity field will only be $\mathcal{C}^0$ across the interface $\Sigma$. This is a consequence of the momentum jump conditions, which usually lead to a jump of $\nabla v$ at the interface. Here we assume $v \in \mathcal{C}^1(\RR \times \Omega)$ for simplicity.}. Moreover, we consider a  $\mathcal{C}^2$-hypersurface $\Sigma_0$ which serves as the initial condition for the advection problem. We assume that $\Sigma_0$ is either a closed hypersurface (without boundary) or that its boundary is contained in the domain boundary $\partial\Omega$ (see Fig.~\ref{fig:notation})
\[ \partial\Sigma_0 \subset \partial\Omega. \]
If the boundary of the interface $\partial\Sigma \subset \partial\Omega$ is non-empty, we call it ``contact line'' and denote it with the symbol $\Gamma$. We  do not go into details of the contact line modeling here and only assume that the velocity field $v$ satisfies an impermeability condition at the boundary, i.e.\
\begin{align}\label{eqn:impermeability_condition}
v \cdot \ndomega = 0 \quad \text{at} \ \partial\Omega, 
\end{align}
where $\ndomega$ denotes the unit outer normal field to $\partial\Omega$.\\
\\
\textbf{Problem statement:} Given $\Sigma_0$ and $v=v(t,x)$ as described above, find the moving hypersurface 
\[ \mathcal{M} = \gr(\Sigma) =  \bigcup_{t \in \RR^+_0} \{t\} \times \Sigma(t) \in \mathcal{C}^1(\RR \times \Omega) \]
that satisfies the initial condition
\begin{align}\label{eqn:initial_condition}
\Sigma(0) = \Sigma_0 
\end{align}
and the kinematic condition
\begin{align}\label{eqn:kinematic_condition}
\normalspeed = v \cdot \nsigma \quad \text{on} \ \gr(\Sigma).
\end{align}
Here, the symbol $\normalspeed$ denotes the speed of normal displacement of $\Sigma$ (see Appendix~\ref{sec:appendix_interfaces_and_kinematics}).\\
\\
Note that the moving hypersurface $\mathcal{M}$ is uniquely determined by \eqref{eqn:initial_condition} and \eqref{eqn:kinematic_condition} if the velocity field $v$ is sufficiently regular and tangential to the boundary as required in \eqref{eqn:impermeability_condition}; see \cite{Fricke2019,Bothe2020,Fricke2021} for more details. In particular, we emphasize there is no additional boundary condition for the interface orientation to be satisfied. For an introduction to the mathematical modeling of dynamic contact lines, the reader is referred to \cite{Bonn2009},\cite{Shikhmurzaev2008},\cite{Fricke2021}.\\
\\
So in the standard level set method, we consider the hyperbolic initial value problem
\begin{align*}
\partial_t \phi + v \cdot \nabla \phi &= 0, \quad x \in \Omega, \ t > 0,\\
\phi(0,x) &= \phi_0(x), \quad x \in \Omega,
\end{align*}
where $\phi_0$ is (locally) a signed distance field of $\Sigma_0$ and $v$ satisfies \eqref{eqn:impermeability_condition}.

\subsection{Preliminary remarks}

\paragraph{Definition of the normal field:} Note that we can define a natural extension of the normal field away from $\Sigma$ by means of
\begin{align}
n = \frac{\nabla \phi}{|\nabla \phi|}
\end{align}
as long as $\nabla \phi \neq 0$. Clearly, the extended field will coincide with $\nsigma$ at the zero contour, i.e.\ $n=\nsigma$ at $\phi=0$.

\paragraph{Rate-of-change of the norm of the gradient:} As already mentioned a few times above, the gradient norm $|\nabla \phi|$ is not conserved by a solution of \eqref{eqn:level_set_equation}. We compute the rate-of-change of $|\nabla \phi|$ explicitly, where we use the simple identity
\begin{align}\label{eqn:gradient_norm_help_identity}
\DDT{} |\nabla \phi| = \DDT{} \left( |\nabla \phi|^2 \right)^{1/2} = \frac{1}{|\nabla \phi|} \inproduct{\nabla \phi}{\DDT{} \nabla \phi} 
\end{align}
Application of the gradient operator to equation~\eqref{eqn:level_set_equation} shows that (see \cite{Fricke2018,Fricke2019})
\begin{align}\label{eqn:ddt_grad_phi}
\DDT{} \nabla \phi = - (\nabla v)^\transpose \nabla \phi.
\end{align}
Hence, it follows that
\begin{equation}\label{eqn:change-of-norm-standard-levelset}
\begin{aligned}
\DDT{} |\nabla \phi| &= - |\nabla \phi| \inproduct{(\nabla v) \frac{\nabla \phi}{|\nabla \phi|}}{\frac{\nabla \phi}{|\nabla \phi|}}\\
&= - |\nabla \phi| \inproduct{(\nabla v) \, n}{n} =: |\nabla \phi| \, \mathcal{R}.
\end{aligned}
\end{equation}

\begin{remark}[Rate of interface generation]
Note that the term
\begin{align}\label{eqn:interface_generation_rate}
\nabla_\Sigma \cdot v = \nabla \cdot v - \inproduct{(\nabla v) \nsigma}{\nsigma}
\end{align}
has an important physical interpretation. In fact, it is nothing but the local rate of \emph{interface generation}. Hence, it measures how much an infinitesimal co-moving control area on the interface is expanded or compressed while it is transported by the flow. This can be seen from the surface transport theorem for co-moving areas (see Appendix~\ref{section:surface_transport}) as follows
$$
\frac{\mathrm{d}}{\mathrm{d} t}|A(t)|=\frac{\mathrm{d}}{\mathrm{d} t} \int_{A(t)} \mathrm{d} o=\int_{A(t)} \operatorname{div}_{\Sigma} \mathbf{v} \mathrm{d} o.
$$
\change{Clearly, in the case of incompressible flows, the interface generation rate $\nabla_\Sigma \cdot v$ coincides with the function $\mathcal{R}$ defined in equation~\eqref{eqn:change-of-norm-standard-levelset}, i.e.}
\begin{align}
\mathcal{R} = - \inproduct{(\nabla v) \nsigma}{\nsigma} = \nabla_\Sigma \cdot v.
\end{align}
Note that equation \eqref{eqn:interface_generation_rate} implies that the quantity $\inproduct{(\nabla v) n}{n}$ must be continuous across the interface in an incompressible two-phase flow setting (since it approaches \change{the interface generation rate} from both sides). Consequently, the function $\mathcal{R}$ has a natural extension away from the interface by means of $n=\nabla \phi/|\nabla \phi|.$ We denote this extension again as $\mathcal{R}$, i.e.\ we let
\[ \mathcal{R} = - \inproduct{(\nabla v) n}{n}(t,x) \quad \text{for} \ \phi \neq 0. \]
\end{remark} 

\subsection{Derivation of the modified level set advection equation}
Motivated by the previous considerations, we derive a transport equation for the level set function that preserves the norm of the gradient at the zero contour. Evidently, the dynamics of the zero contour must still be the same as in the standard formulation. To achieve this, we consider the following ansatz
\begin{align}\label{eqn:levelset-source-term-generic}
\partial_t \phistar + v \cdot \nabla \phistar = \phistar f ( \phi ( \cdot), v(\cdot )).
\end{align}
Equation \eqref{eqn:levelset-source-term-generic} leaves the motion of the zero contour invariant provided that $f$ is bounded as $\phi \rightarrow 0$. On the other hand, the function $f$ can be carefully chosen in order to achieve
\begin{align}
\DDT{} |\nabla \phistar|\Big|_{\phistar=0} = 0.
\end{align}

For this purpose, we compute the rate-of-change of $|\nabla \phistar|$ for a solution of \eqref{eqn:levelset-source-term-generic}. Application of the gradient operator to a solution of \eqref{eqn:levelset-source-term-generic} yields
\begin{align*}
\DDT{} \nabla \phistar = - (\nabla v)^\transpose \nabla \phistar + f(\phistar,v) \nabla \phistar + \phistar \, \nabla f(\phistar,v).
\end{align*}
Hence, we obtain using \eqref{eqn:gradient_norm_help_identity} and $n=\nabla \phistar /|\nabla \phistar|$
\begin{align}\label{eqn:change-of-norm-generic}
\DDT{} |\nabla \phistar| = |\nabla \phistar| \left( - \inproduct{(\nabla v)n}{n} + f\right) + \phistar n \cdot \nabla f.
\end{align}
Now, we have to choose $f$ in such a way that the right-hand side of \eqref{eqn:change-of-norm-generic} vanishes for $\phistar=0$. Obviously, the choice
\begin{align}
f(\phistar,v) := - \mathcal{R} = \inproduct{(\nabla v)n}{n}
\end{align}
does the job. In this case, we have
\begin{equation}
\begin{aligned}
\DDT{} |\nabla \phistar| = - \phistar \frac{\partial \mathcal{R}}{\partial n}.
\end{aligned}
\end{equation}
In particular, the norm of the gradient of $\phistar$ is preserved at the interface. 

\paragraph{Summary:} We propose to investigate the non-linear equation\footnote{One of us (IR) derived and presented this modified level set equation already in 2011 as part of a scientific workshop \cite{Roisman2011}. Within the present study, it has been derived independently (by MF) in the form as presented here. Recently, Hamamuki presented an analytical investigation on the modified equation in the framework of viscosity solutions \cite{Hamamuki2019}. We are grateful to Prof.\ Y.\ Giga for pointing out this reference.}
\begin{equation}\label{eqn:adapted_levelset_equation} 
\begin{aligned}
\partial_t \phistar + v \cdot \nabla \phistar &= \phistar \inproduct{(\nabla v) \frac{\nabla \phistar}{|\nabla \phistar|}}{\frac{\nabla \phistar}{|\nabla \phistar|}} = - \phistar \mathcal{R}.
\end{aligned}
\end{equation}

\section{An adapted upwind method}\label{sec:upwind}
In order to construct numerical methods based on \eqref{eqn:adapted_levelset_equation}, we note that a reasonable approach is to approximate the interface generation rate $\mathcal{R}$ in each time interval $[t^n,t^{n+1}]$ by some function $r=r(t,x)$ that is decoupled from the level set function $\phi$. For example, $r$ could be a snapshot of $\mathcal{R}$ at time $t^n$. Thanks to this approach, we will consider the \emph{linear} problem
\begin{align}\label{eqn:adapted_levelset_equation_linearized}
\partial_t \phi + v \cdot \nabla \phi = - r(x) \, \phi, \quad x \in \Omega, \ t > 0
\end{align}
in each time step. In fact, the zero-contour of the solution of the linear problem \eqref{eqn:adapted_levelset_equation_linearized} coincides with the zero-contour of the non-linear problem \eqref{eqn:adapted_levelset_equation}. Hence, approximating the interface generation rate $\mathcal{R}(t,x)$ by $r(x)$ only leads to deviations in $|\nabla \phi|$ with no impact on the evolution of the zero-contour. Hence, on the continuous level, no accuracy is lost in the position of the interface. On the other hand, small deviations from $|\nabla \phi|=1$ at the interface are usually unproblematic in practice. Moreover, the normalization of the gradient at the interface is consistently achieved when refining the mesh and the timestep. Therefore, we may concentrate on linear problems of type \eqref{eqn:adapted_levelset_equation_linearized} for the construction of numerical methods.

\begin{remark}
Using the method of characteristics (see, e.g., \cite{Evans2010}), it is easy to show that equation \eqref{eqn:adapted_levelset_equation_linearized} is well-posed as an initial value problem if the velocity field satisfies the impermeability condition
\[ v \cdot \ndomega = 0 \quad \text{at} \ \partial\Omega. \]
The reason is that the characteristics of \eqref{eqn:level_set_equation} and \eqref{eqn:adapted_levelset_equation_linearized} are the same. The only difference is that for \eqref{eqn:adapted_levelset_equation_linearized}, the function $\phi$ is no longer constant along characteristics. Hence, no boundary condition has to be imposed for \eqref{eqn:adapted_levelset_equation_linearized} and the kinematics of the contact angle is not affected. This makes this numerical approach particularly interesting for the simulation of dynamic wetting problems.
\end{remark}

\paragraph{Derivation of the numerical scheme:} Integrating \eqref{eqn:adapted_levelset_equation_linearized} over a control volume $V_{ijk}$ yields (assuming $\nabla \cdot v = 0$)
\begin{align}\label{eqn:adapted_levelset_equation_linearized_integral}
\ddt{} \frac{1}{|V_{ijk}|} \int_{V_{ijk}} \phi \, dV = \frac{1}{|V_{ijk}|} \int_{\partial V_{ijk}} \phi \, v \cdot n \, dA - \frac{1}{|V_{ijk}|} \int_{V_{ijk}} r \phi \, dV. 
\end{align}
One of the simplest possible methods to solve \eqref{eqn:adapted_levelset_equation_linearized_integral} is the first-order upwind scheme. It can be seen as a primer for more sophisticated methods with a higher order of accuracy. In order to derive it, we may approximate the last term as
\[ \frac{1}{|V_{ijk}|} \int_{V_{ijk}} r \phi \, dV \approx \left(\frac{1}{|V_{ijk}|} \int_{V_{ijk}} r \, dV \right) \left( \frac{1}{|V_{ijk}|} \int_{V_{ijk}} \phi \, dV \right).   \]
With this approximation, we obtain a particularly well-suited problem to be solved with finite volume methods
\begin{align}\label{eqn:adapted_levelset_equation_linearized_integral_v2}
\ddt{} \frac{1}{|V_{ijk}|} \int_{V_{ijk}} \phi \, dV = \frac{1}{|V_{ijk}|} \int_{\partial V_{ijk}} \phi \, v \cdot n \, dA - \left(\frac{1}{|V_{ijk}|} \int_{V_{ijk}} r \, dV \right) \left( \frac{1}{|V_{ijk}|} \int_{V_{ijk}} \phi \, dV \right).
\end{align}
A simple time-explicit semi-discretization of \eqref{eqn:adapted_levelset_equation_linearized_integral_v2} reads as
\begin{equation}\label{eqn:adapted_first_order_upwind}
\begin{aligned}
\frac{\phi_{ijk}^{n+1} - \phi_{ijk}^n}{\Delta t} = \mathcal{F}_{ijk}^n - r_{ijk}^n \, \phi_{ijk}^n \\
\quad \Leftrightarrow \quad \phi_{ijk}^{n+1} = \phi_{ijk}^{n} (1 - r_{ijk}^n \Delta t) + \Delta t \, \mathcal{F}_{ijk}^n.
\end{aligned}
\end{equation}
where $\mathcal{F}_{ijk}^n$ denotes the total (sum over all faces) upwind flux\footnote{See \cite{LeVeque2002} for a discussion of the first-order upwind method.} for the control volume $V_{ijk}$ at time $t_n$ and $r_{ijk}^n$ is a discretization of the rate of interface generation. From equation \eqref{eqn:adapted_first_order_upwind}, we see that the numerical time step should be chosen such that $|r_{ijk}^n| \Delta t \leq C_r < 1$. This ensures that the time scale associated with the source term is resolved.\newline
\newline
Notably, the only implementation work needed to extend the standard upwind scheme for \eqref{eqn:adapted_levelset_equation} is to discretize the source term
\[ r_{ijk}^n = -\inproduct{(\nabla v) \, n}{n}_{ijk}^n = - \inproduct{(\nabla v) \, \frac{\nabla \phi}{|\nabla \phi|}}{\frac{\nabla \phi}{|\nabla \phi|}}_{ijk}^n \]
and to include it in the time integration according to \eqref{eqn:adapted_first_order_upwind}. For the purpose of this demonstrator method, we take the analytical velocity gradient at the cell center and discretize $\nabla \phi$ with second-order central finite differences away from the boundary. At boundaries, we use the second-order forward-/backward finite difference schemes based on the Taylor expansions
\begin{equation}
\begin{aligned}
f'(x) &= \frac{-3 f(x) + 4 f(x + \Delta x) - f(x+2 \Delta x)}{2 \Delta x} + \mathcal{O}(\Delta x^2),\\
f'(x) &= \frac{3 f(x) - 4 f(x - \Delta x) + f(x-2 \Delta x)}{2 \Delta x} + \mathcal{O}(\Delta x^2).
\end{aligned}
\end{equation}
Second-order accuracy is reached if the function $f$ is of class $\mathcal{C}^3$. 

\begin{remark}[Regularization of the source term] 
\change{There are two important remarks regarding the numerical treatment of the source term.}
\begin{enumerate}[(i)]
 \item Obviously, the source term $r_{ijk}^n$ is not well-defined if $|\nabla \phi| = 0$. Therefore, it must be regularized in these regions. To achieve this, we replace
\[ \frac{\nabla \phi}{|\nabla \phi|} \quad \text{by} \quad \frac{\nabla \phi}{|\nabla \phi|+\varepsilon} \]
for some small regularization parameter $\varepsilon$. This will affect the calculation of the normal vector only in regions where $|\nabla \phi|$ is already very small. We use $\varepsilon = 10^{-12}$ for our numerical experiments in Section~\ref{sec:results}.
\item \change{Moreover, there may be a problem at inflow boundaries introduced by the source term. A discontinuity in the normal field is generated at inflow boundaries if the homogeneous Neumann condition is applied for the level set field. This would in turn degenerate the regularity of the function $r(x)$ and create numerical instabilities. To avoid this problem, we use a numerical cu-toff (or Mollifier) function which forces the source term to zero sufficiently far away from the zero contour. So, in summary, we set}
\begin{align}
r_{ijk}^n = - \inproduct{(\nabla v) \, \frac{\nabla \phi}{|\nabla \phi|+\varepsilon}}{\frac{\nabla \phi}{|\nabla \phi|+\varepsilon}}_{ijk}^n G(\phi_{ijk}^n),
\end{align}
\change{where $G$ is the symmetric, $\mathcal{C}^1$-regular cut-off function defined as}
\begin{align}\label{eqn:mollifier_function}
G(x) = \begin{cases} 1 & \text{if} \quad 0 \leq x \leq w_1, \\
\exp\left(- \ln(10^3) \frac{(x-w_1)^2}{(w_2-w_1)^2}\right) & \text{if} \quad w_1 < x, \\
G(-x) & \text{if} \quad x < 0.
\end{cases}
\end{align}
\change{Here $w_1$ denotes the width of the region where $G$ is unity and $w_2>w_1$ is the distance at which $G=10^{-3}$. It is evident from Figure~\ref{fig:mollifier} that the use of the Mollifier is very important in practice if inflow boundaries are present. In the reported case, all boundaries except for the bottom one act as inflow boundaries (here $v$ is periodic in time).}
\end{enumerate}
\begin{figure}[hb]
\centering
\subfigure[Mollifier inactive.]{\includegraphics[width=0.5\columnwidth]{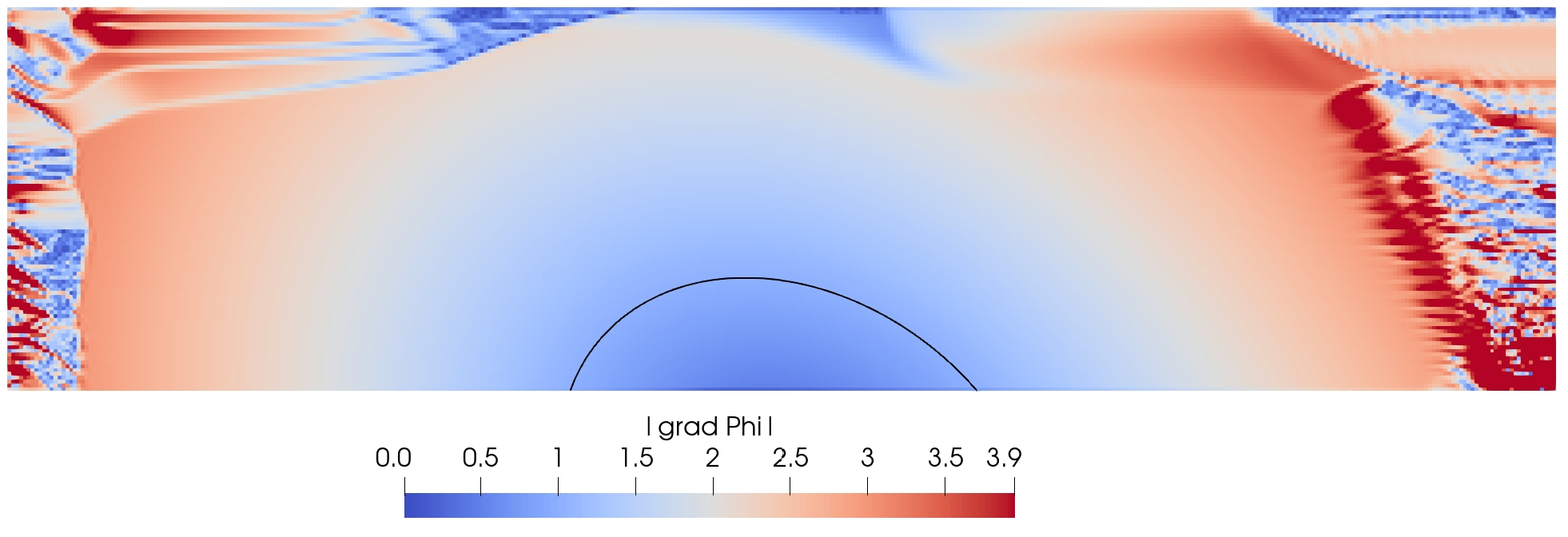}}
\subfigure[Mollifier active.]{\includegraphics[width=0.5\columnwidth]{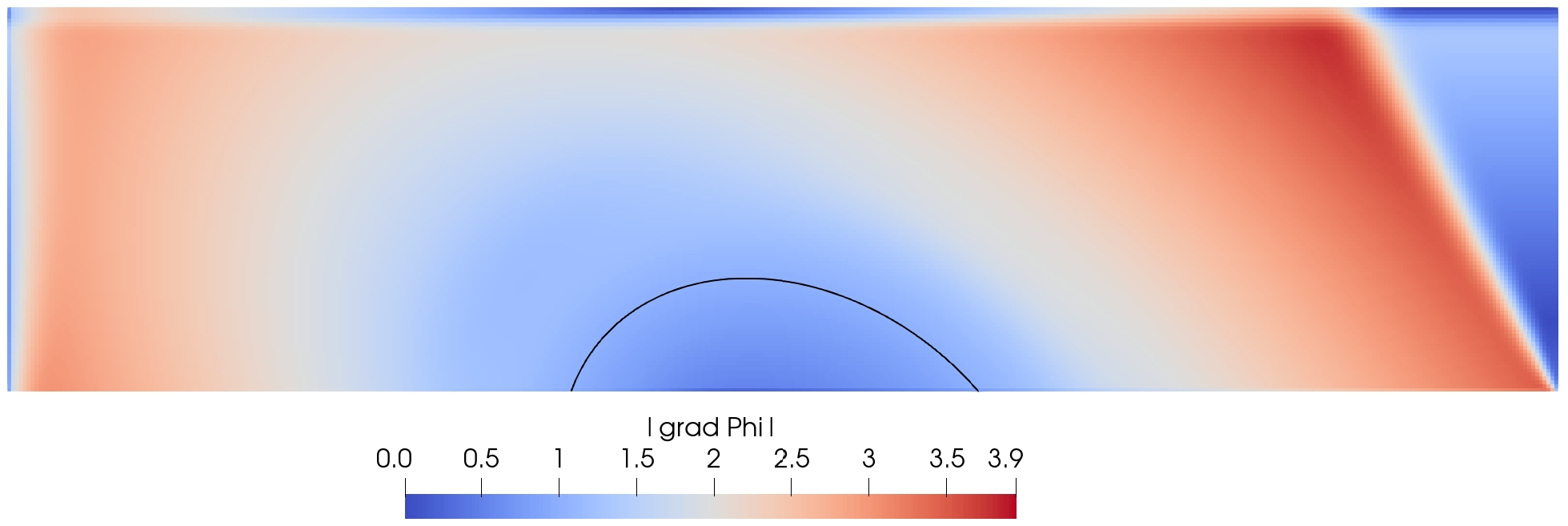}} 
\caption{Graphical comparison of the field $|\nabla \phi|$ with and without mollifier function (inflow/outflow at left, right and top boundary).}\label{fig:mollifier}
\end{figure}
\end{remark}
\vspace{3mm}

The method described above is implemented in a C++ demonstrator code that has been published before (without the source term) as Open Source in \cite{Fricke2019a,Fricke2019a-code}. Note that the method works naturally in both two and three spatial dimensions as the only adaption of the upwind scheme is to discretize the source term.

\clearpage
\section{Numerical examples}\label{sec:results}
\subsection{Numerical examples in two dimensions}
\begin{figure}[ht]
\centering
\includegraphics[width=0.6\columnwidth]{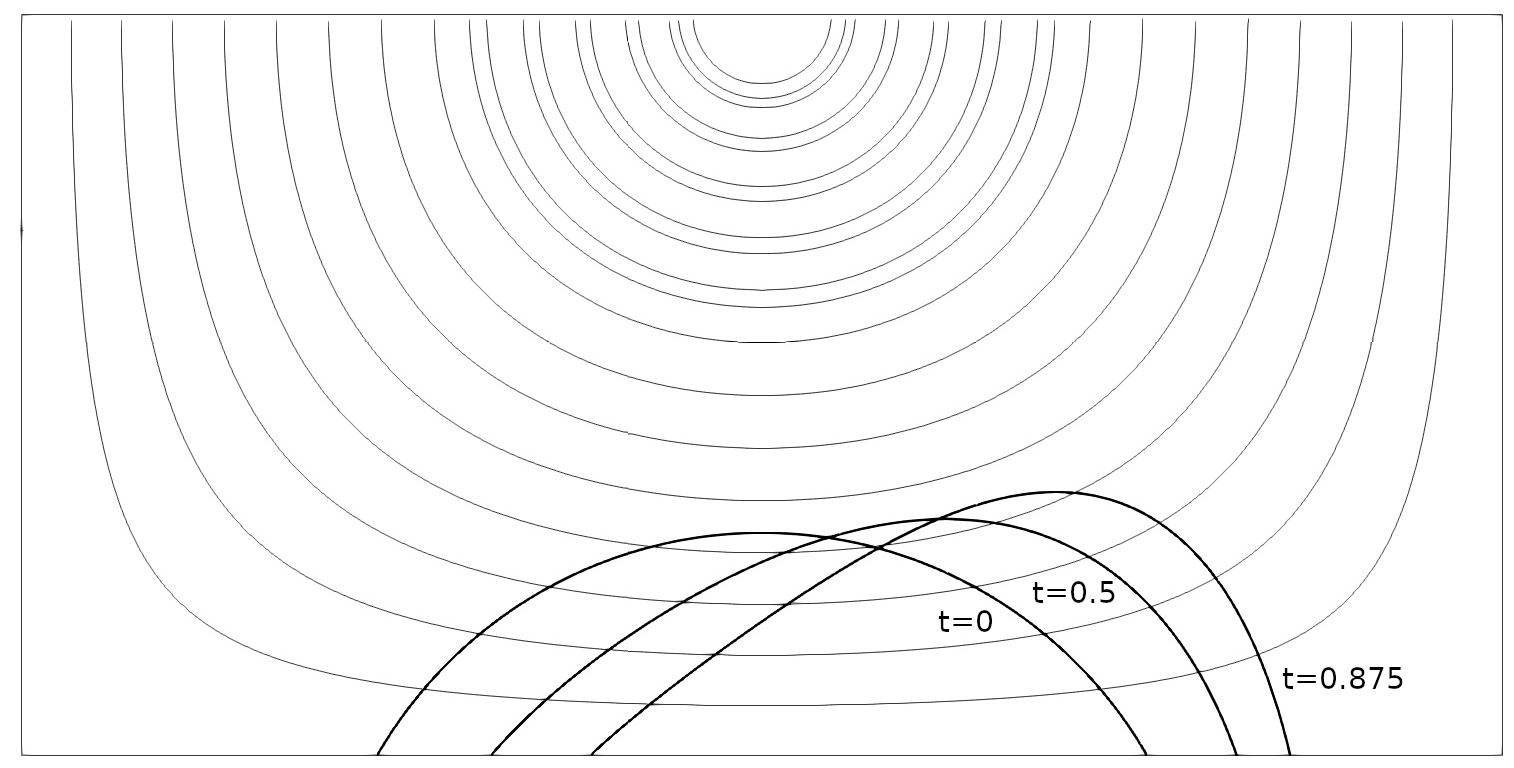}
\caption{Streamlines and snapshots of the interface at $t=0,0.5,0.875$ for the velocity field \eqref{eqn:shear_field}.}\label{fig:dynamics-shear-2d}
\end{figure}
We consider a two-dimensional domain
\[ 0 \leq x \leq 1, \quad 0 \leq y \leq 0.5 \]
and choose the initial level set function
\[ \change{\phi_0(x,y) = \sqrt{(x-x_0)^2 + (y-y_0)^2} - R_0.} \]
for $x_0 = 0.5$, $y_0=-0.15$ and $R_0=0.3$. Hence, the initial zero contour is a semi-circle that meets the boundary $y=0$ at a contact angle of $60$ degrees. Note that $\phi_0$ \change{is precisely the signed distance function, i.e.\ $|\nabla \phi_0|=1$ holds.} For the spatial discretization, we use a uniform Cartesian grid with $\Delta x = \Delta y$. We choose between $100$ and $400$ cells in the $x$-direction corresponding to between $30$ and $120$ cells per initial radius $R_0$. \change{The mollifier function \eqref{eqn:mollifier_function} is used with parameters $w_1=0.05$ and $w_2=0.15$.} \\
\\
Two velocity fields are chosen to transport the level set with the first-order upwind scheme
\begin{enumerate}[(i)]
 \item ``Vortex-in-a-box field'': 
 \begin{align}\label{eqn:shear_field}
  v(x,y) = v_0 (-\sin(\pi x)\cos(\pi y),\cos(\pi x) \sin(\pi y))
 \end{align}
for $v_0 = -0.2$.
\item ``Time-periodic field'':
 \begin{align}\label{eqn:periodic_field}
 v(t,x,y) = \cos\left( \frac{\pi t}{\tau} \right) \, (v_0 + c_1 x + c_2 y, -c_1 y)  
 \end{align}
 for $v_0 = -0.2$, $c_1=0.1$, $c_2=-2$ and $\tau=0.4$.
\end{enumerate}
\begin{figure}[ht]
\subfigure[Source on. Contours for $\phi=0$ and $\phi = \pm w_1$.]{\includegraphics[width=0.5\columnwidth]{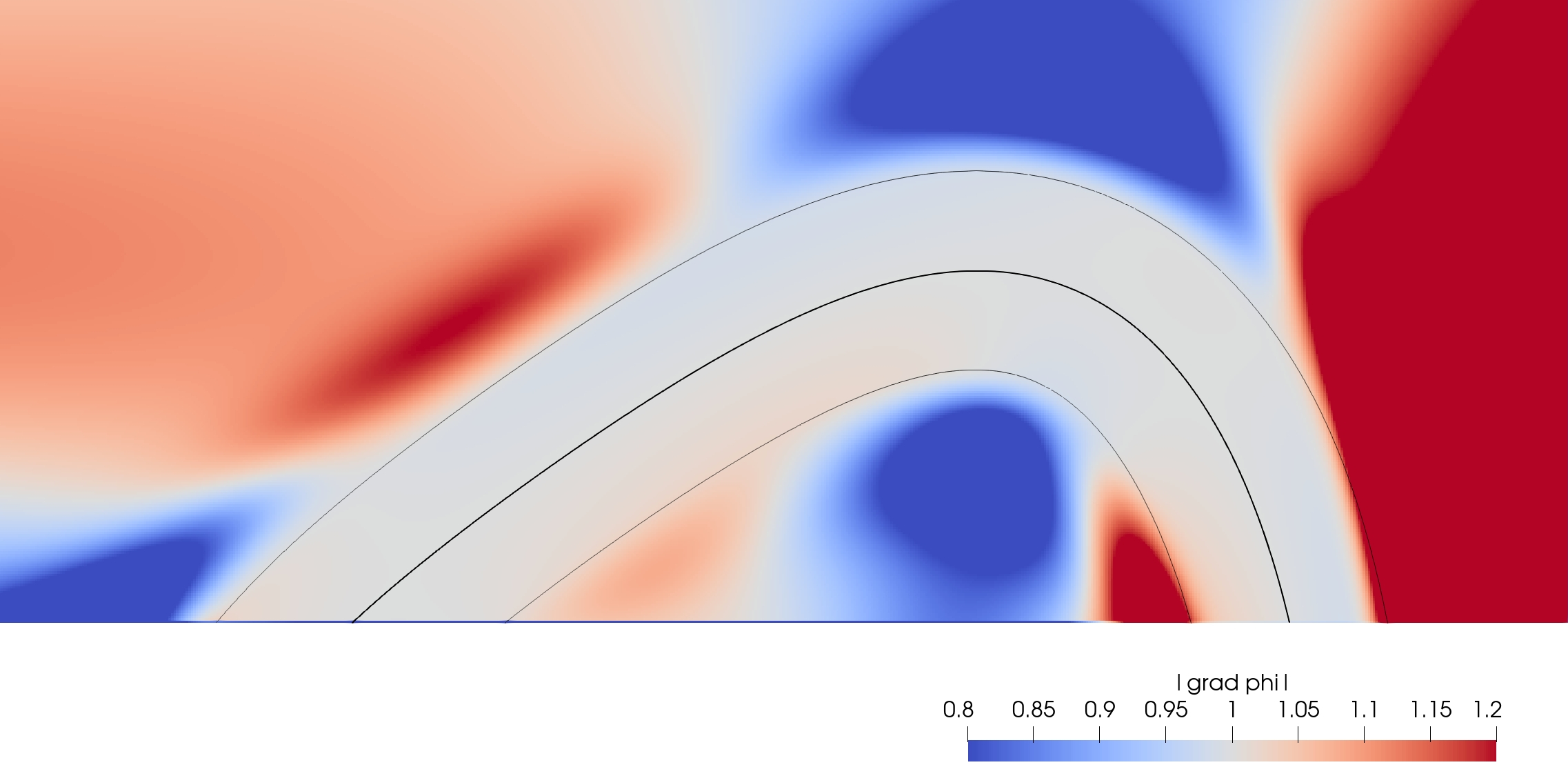}}
\subfigure[Source off. Contour for $\phi=0$.]{\includegraphics[width=0.5\columnwidth]{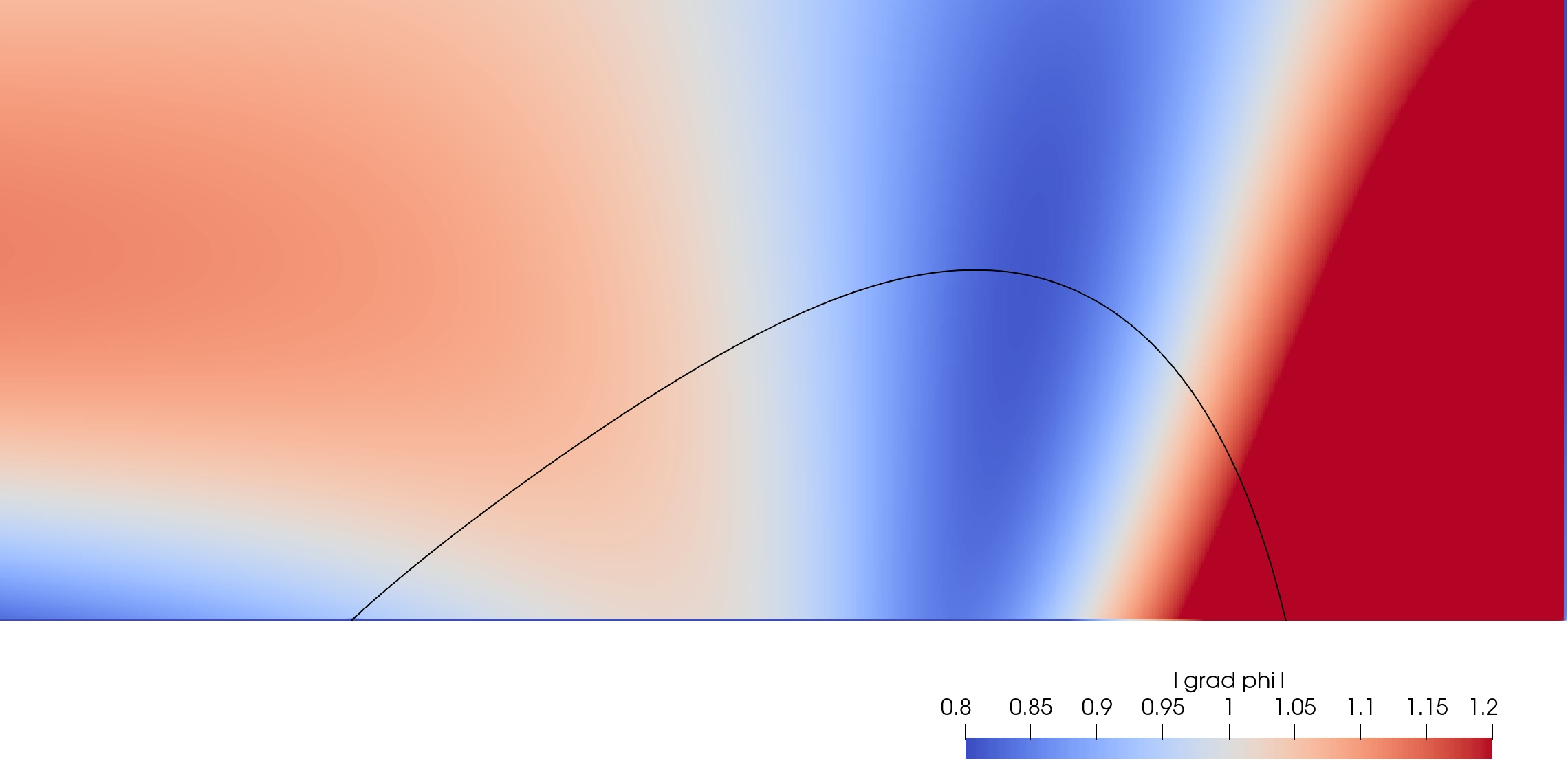}}
\caption{Visualization of the field $|\nabla \phi|$ near the zero contour for the vortex-in-a-box field \eqref{eqn:shear_field} at time $t=0.875$.}\label{fig:norm_grad_phi_field_plot}
\end{figure}

The mesh size is varied to study the numerical convergence of the method. The numerical timestep is varied accordingly such that the CFL number is kept constant ($\text{CFL}=0.5$ in this case). We study the evolution of the ``right'' contact point, initially being located at coordinate $x = x_0 + R \sin \theta_0 = 0.5 + 0.3 \sin(60^\circ) \approx 0.76$. We follow this point along the $y=0$ boundary and evaluate
\begin{enumerate}[(i)]
 \item the coordinate $x(t)$,
 \item the contact angle $\theta(t)$,
 \item the curvature $\kappa(t)$ and
 \item the norm of the gradient $|\nabla \phi|(t)$ for the standard level set method and the new version with active source term.
\end{enumerate}
For the transport of the coordinate $x(t)$, the contact angle $\theta(t)$ as well as the curvature $\kappa(t)$, we have access to reference solutions derived in \cite{Fricke2019,Fricke2020}. In particular, we use the kinematic evolution equation for the dynamic contact angle \cite{Fricke2019} and the kinematic evolution equation for the (mean) curvature (see \cite{Fricke2021}) to study the accuracy of the method.

\begin{figure}[hb]
\subfigure[Gradient norm $|\nabla \phi|(t)$ ($\Delta x = 1.25 \cdot 10^{-3}$).]{\includegraphics[width=0.5\columnwidth]{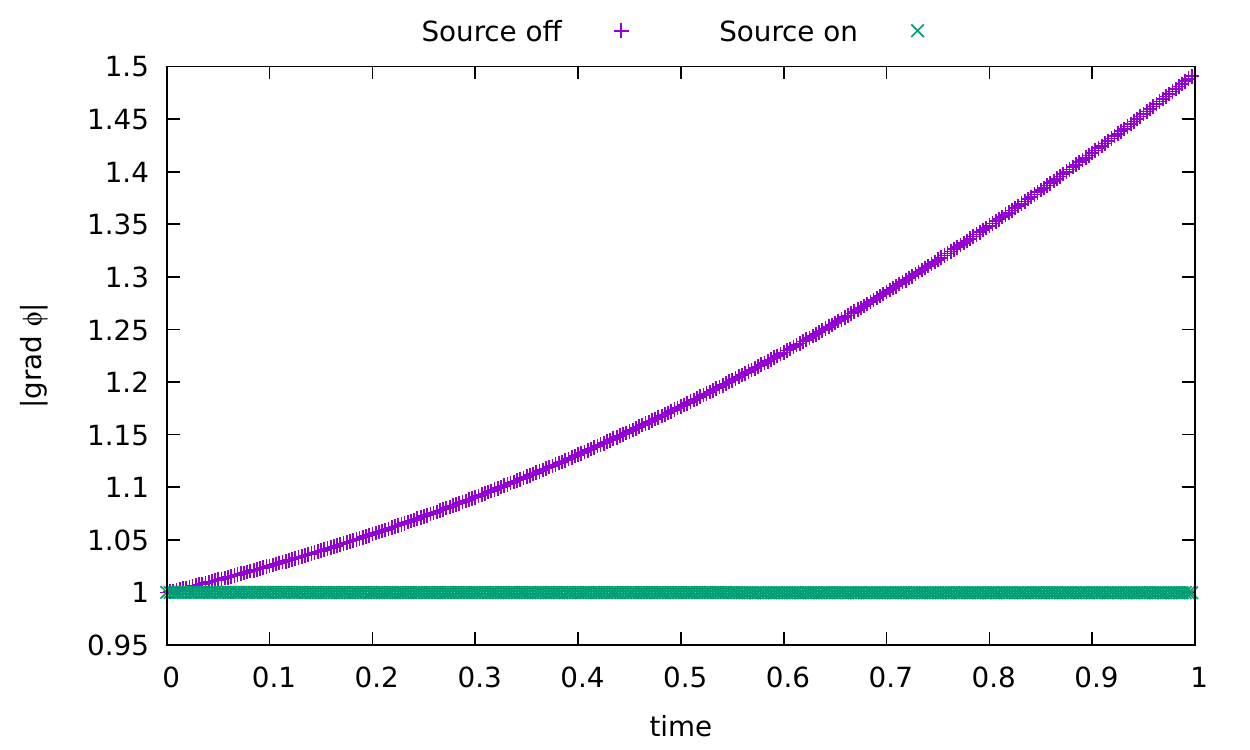}}
\subfigure[Deviation from local signed-distance.]{\includegraphics[width=0.5\columnwidth]{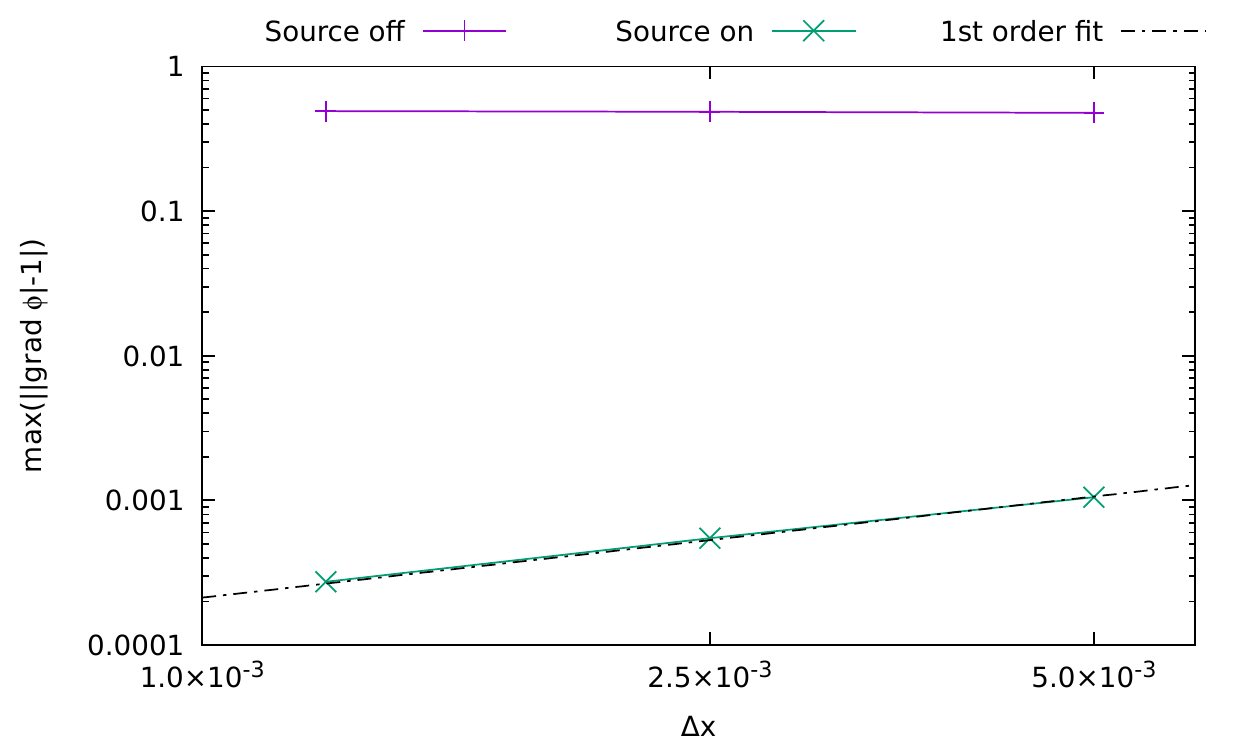}}
\caption{Numerical transport of the gradient norm for the field \eqref{eqn:shear_field}.}
\label{fig:shear/gradient_norm_evolution}
\end{figure}

\begin{figure}[ht]
\subfigure[Position $x(t)$ ($\Delta x = 1.25 \cdot 10^{-3}$).]{\includegraphics[width=0.5\columnwidth]{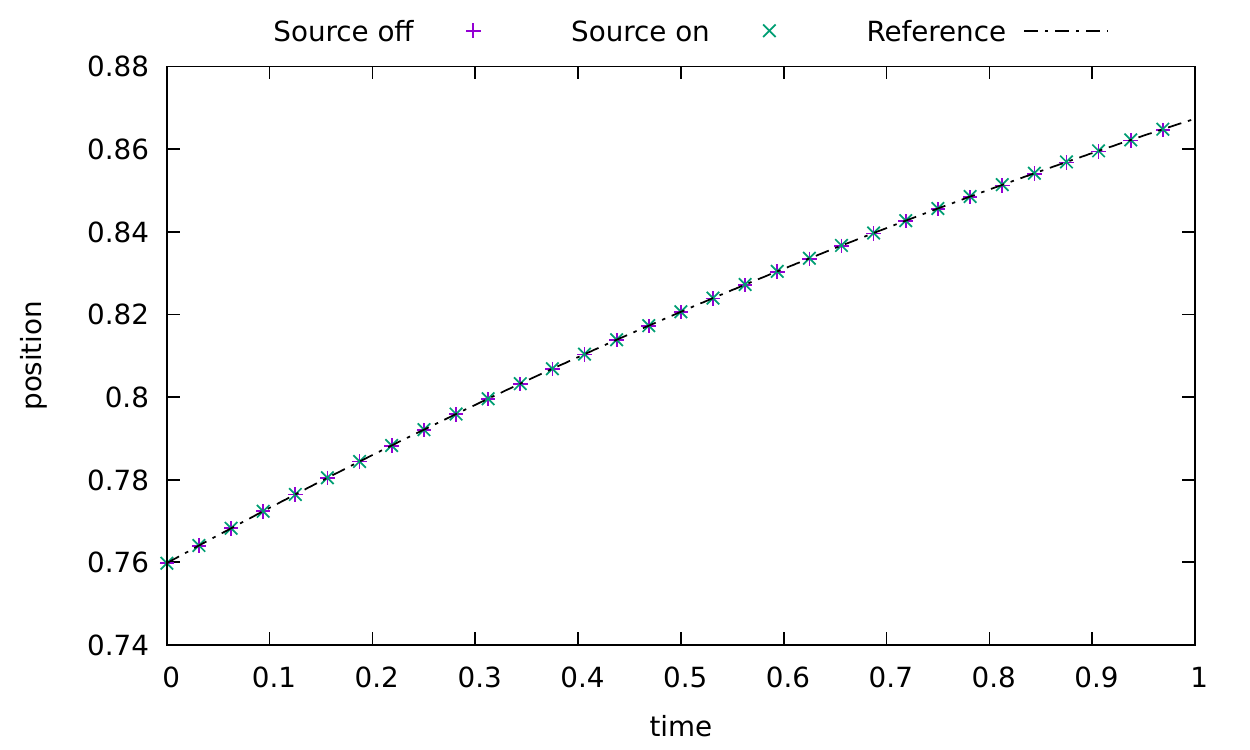}} 
\subfigure[Mesh study for $x(t)$.]{\includegraphics[width=0.5\columnwidth]{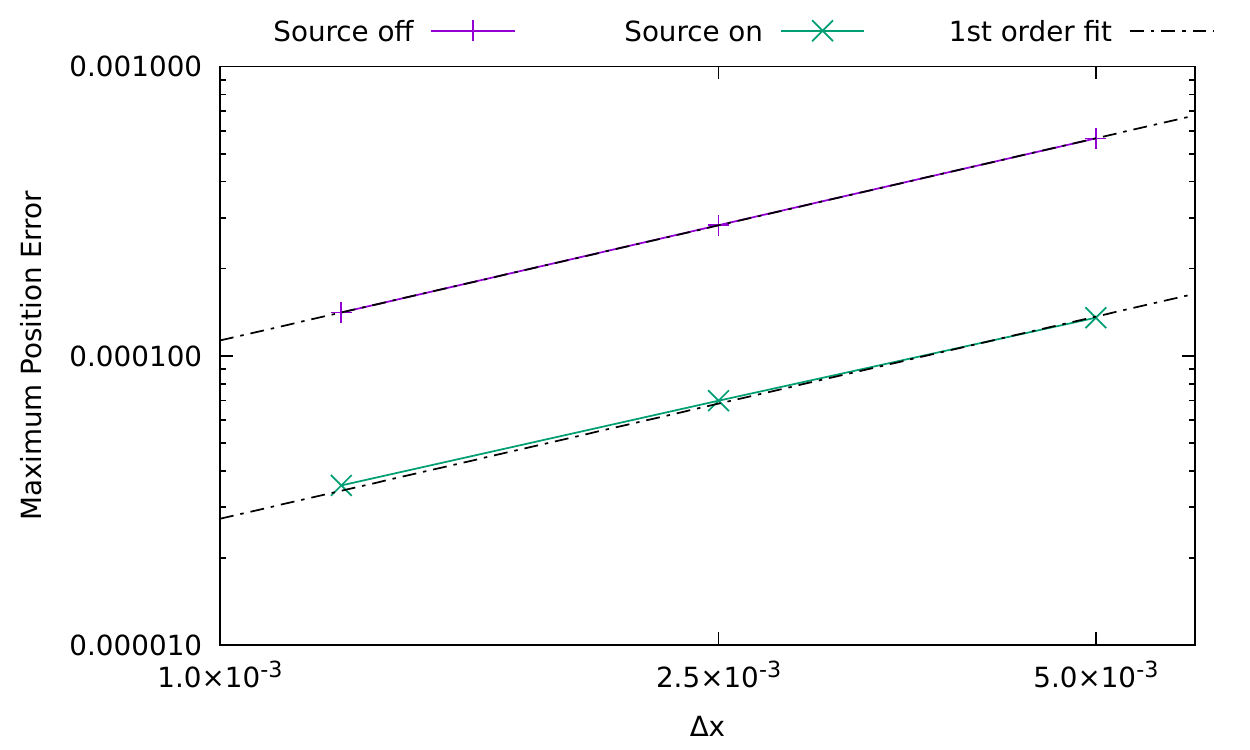}} 
\subfigure[Contact angle $\theta(t)$ ($\Delta x = 1.25 \cdot 10^{-3}$).]{\includegraphics[width=0.5\columnwidth]{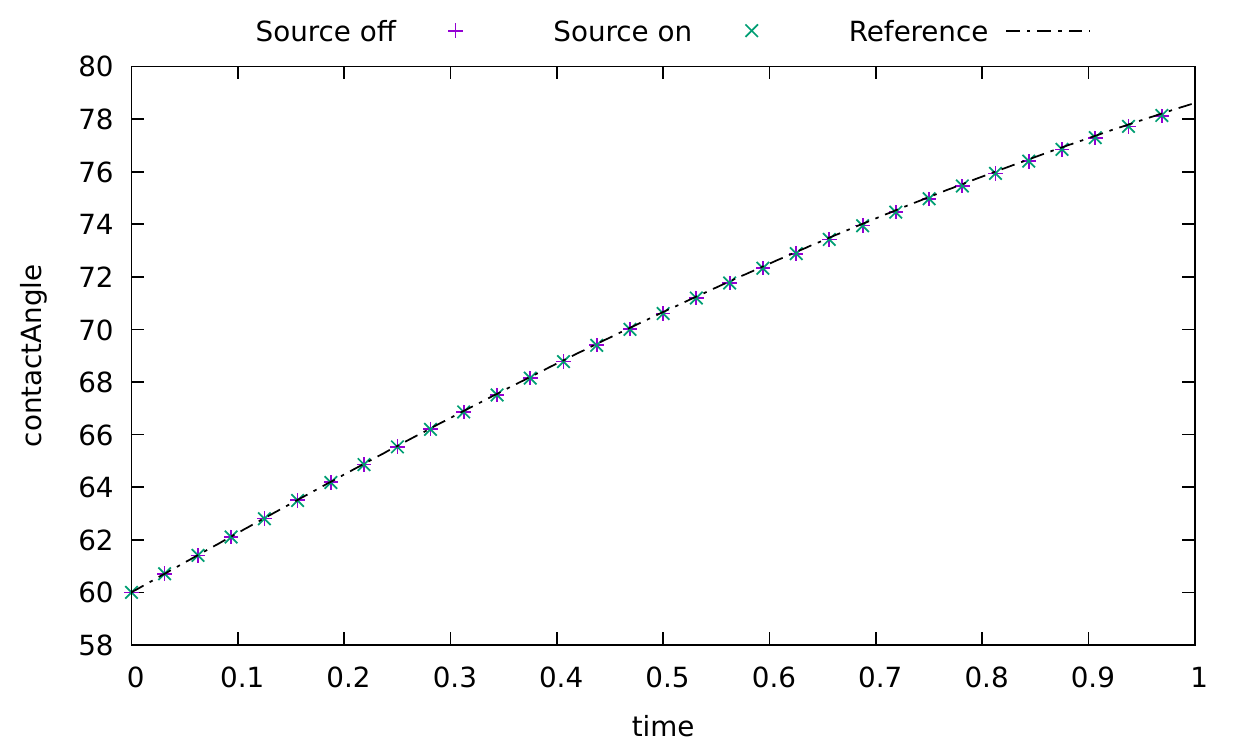}}
\subfigure[Mesh study for $\theta(t)$.]{\includegraphics[width=0.5\columnwidth]{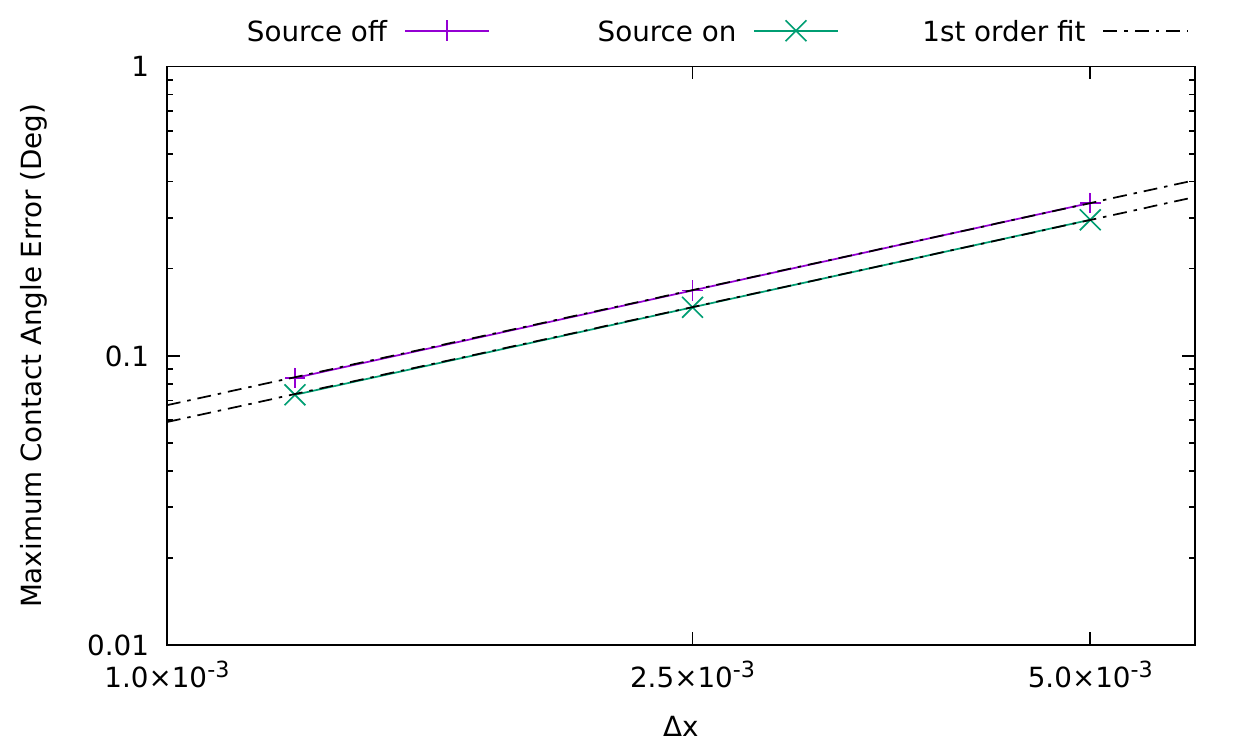}}
\subfigure[Curvature $\kappa(t)$ ($\Delta x = 1.25 \cdot 10^{-3}$).]{\includegraphics[width=0.5\columnwidth]{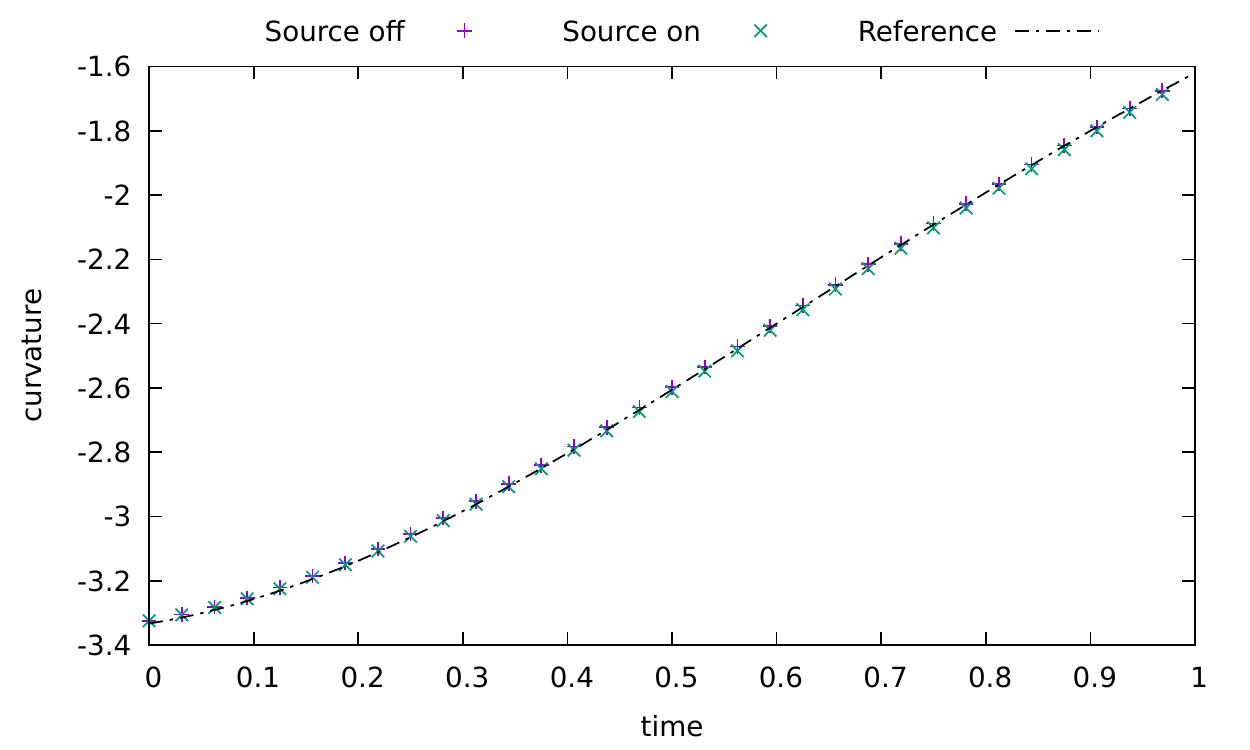}}
\subfigure[Mesh study for $\kappa(t)$.]{\includegraphics[width=0.5\columnwidth]{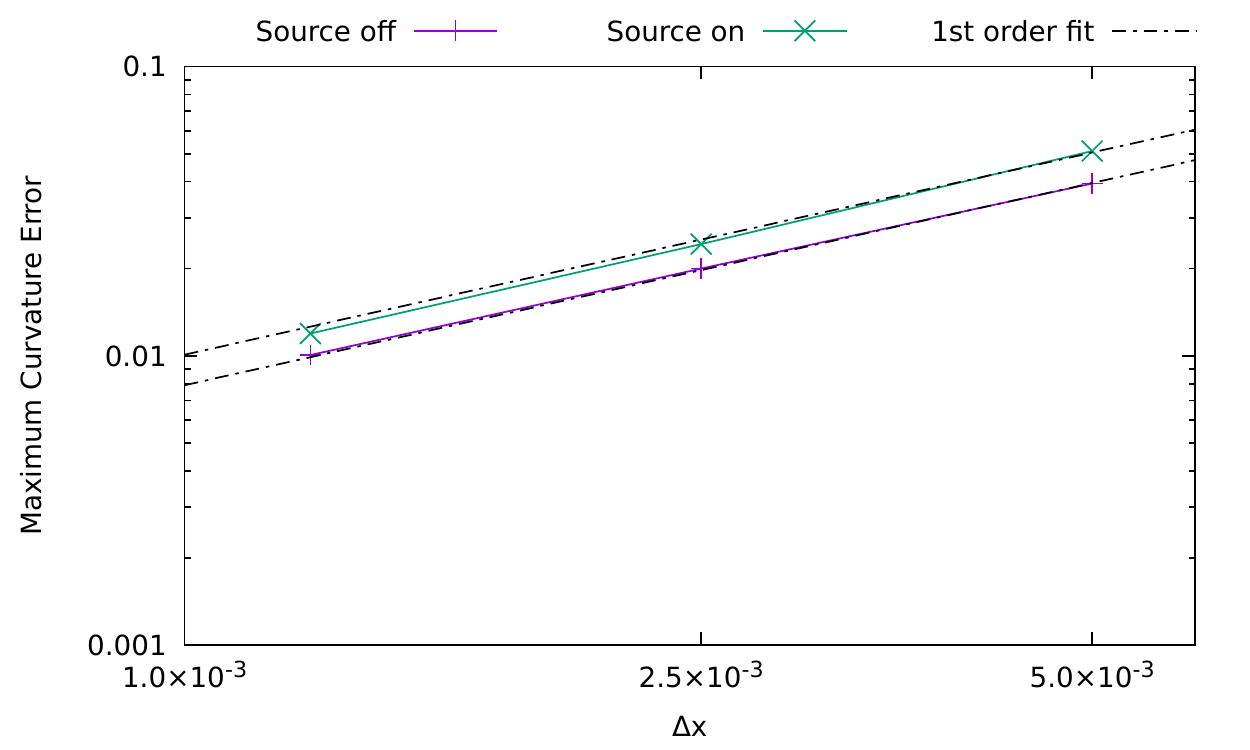}}
\caption{Convergence study for the contact point position, contact angle and curvature for the vortex-in-a-box field \eqref{eqn:shear_field}.}
\label{fig:shear/convergence_study}
\end{figure}

\begin{figure}[ht]
\subfigure[Contact angle $\theta(t)$ ($\Delta x = 1.25 \cdot 10^{-3}$).]{\includegraphics[width=0.5\columnwidth]{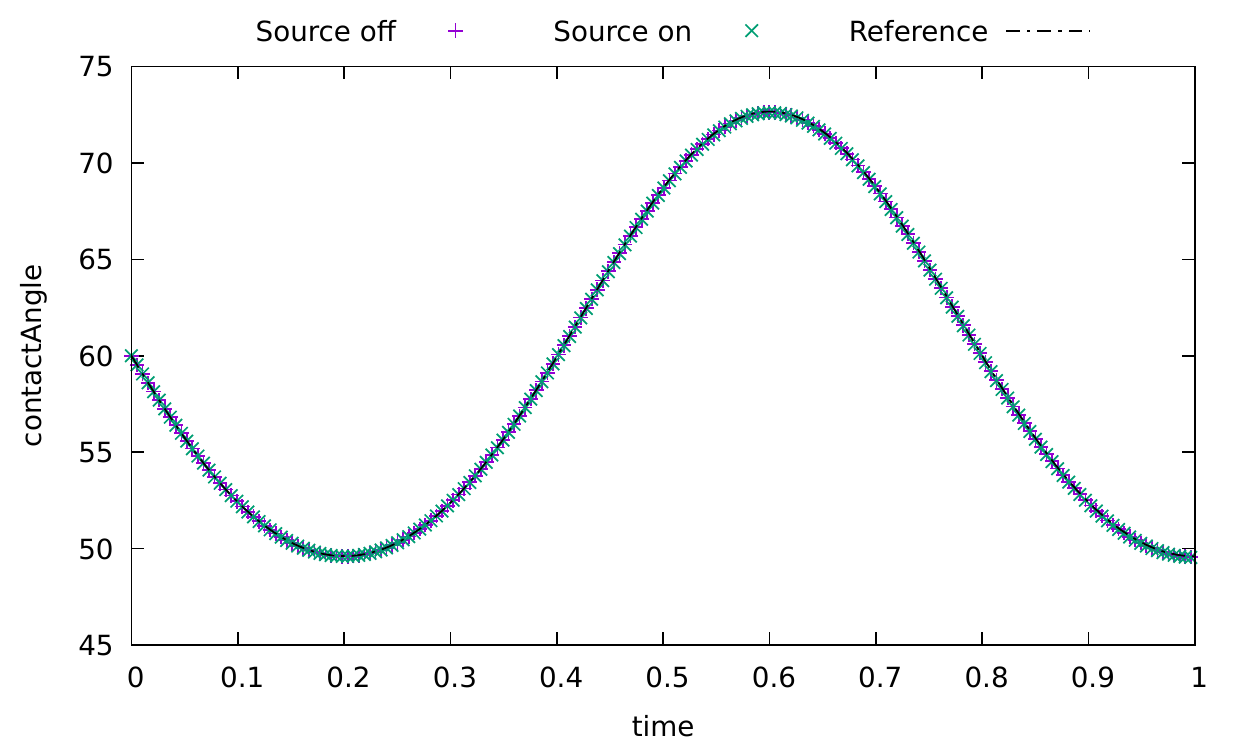}} 
\subfigure[Mesh study for $\theta(t)$.]{\includegraphics[width=0.5\columnwidth]{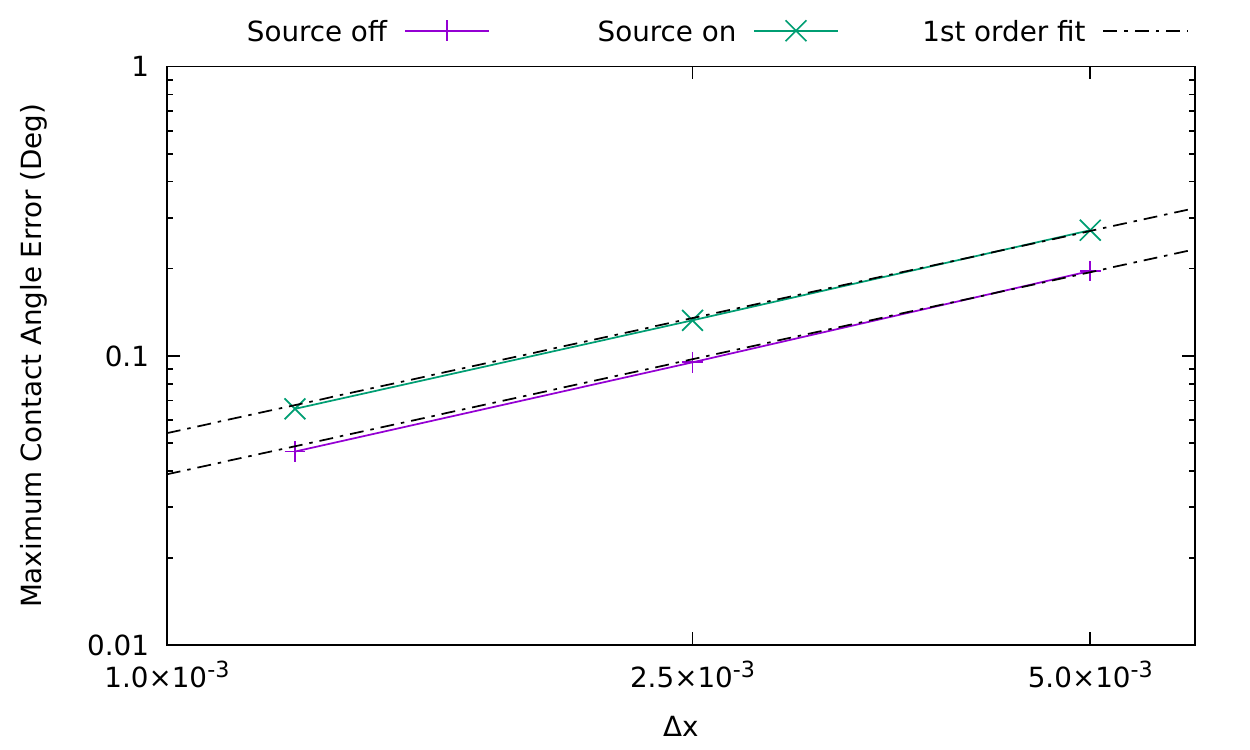}}
\subfigure[Curvature $\kappa(t)$ ($\Delta x = 1.25 \cdot 10^{-3}$).]{\includegraphics[width=0.5\columnwidth]{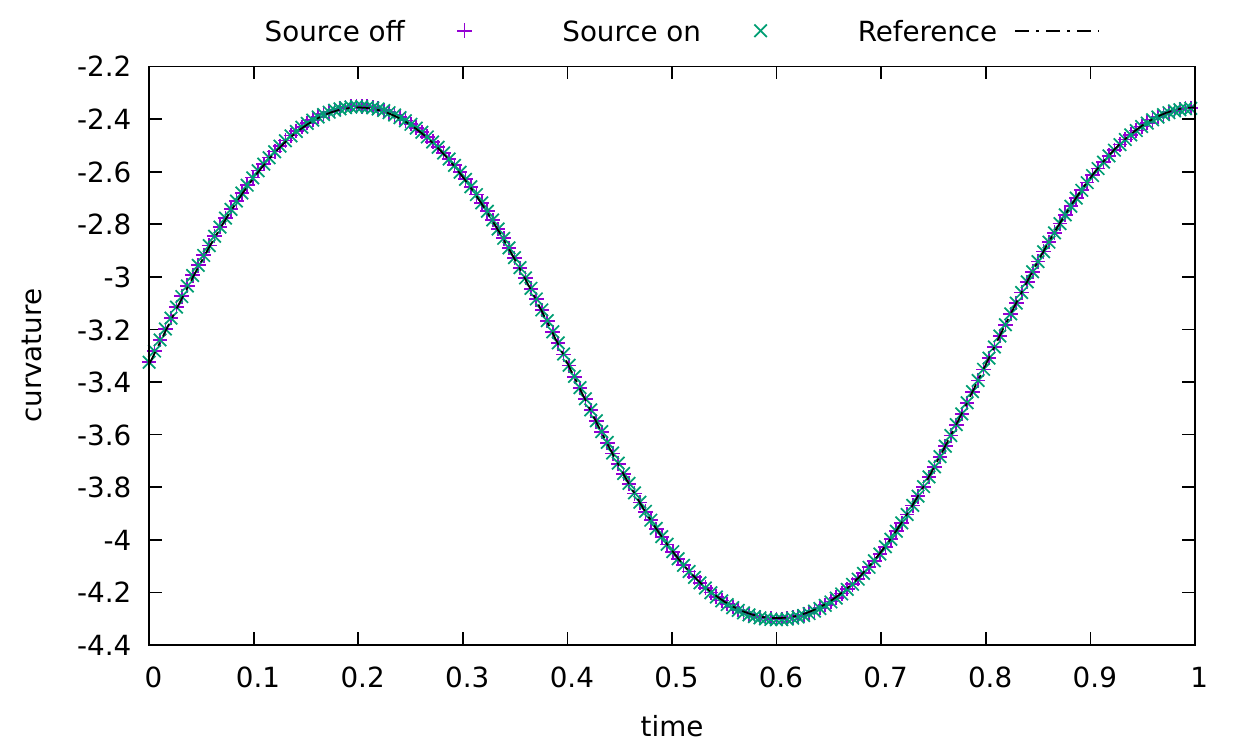}} 
\subfigure[Mesh study for $\kappa(t)$.]{\includegraphics[width=0.5\columnwidth]{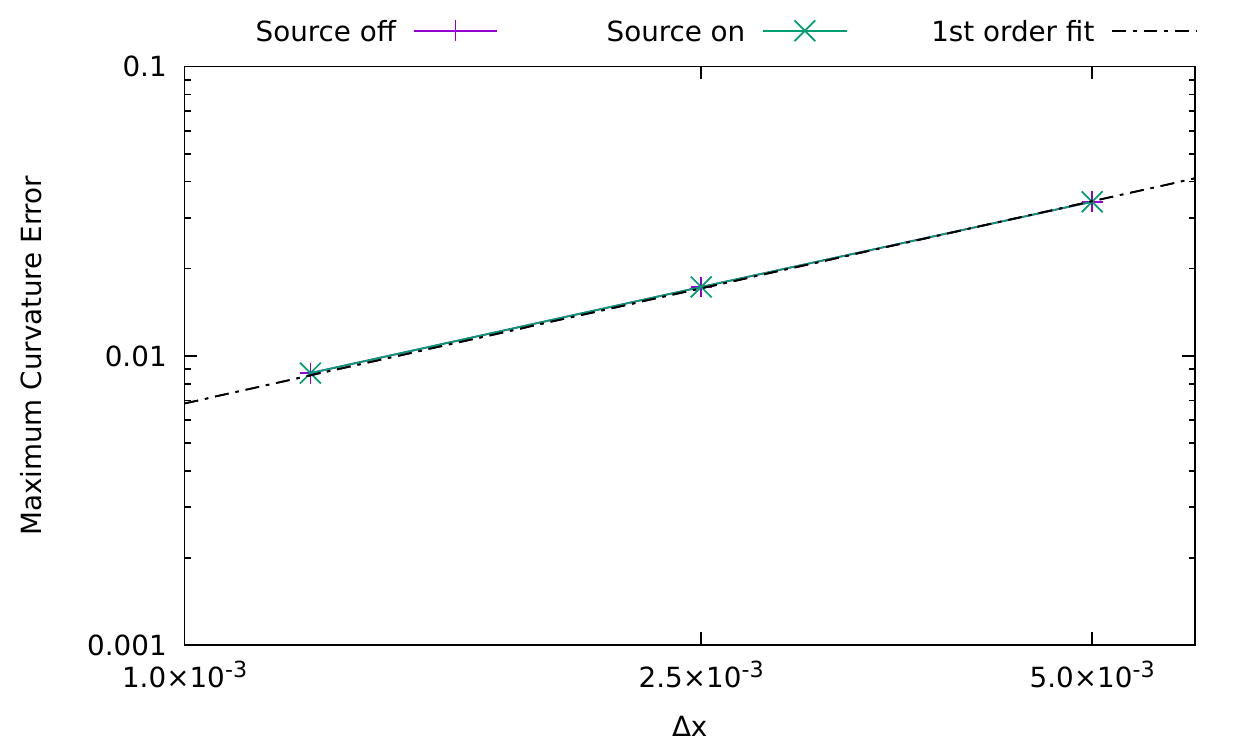}}
\subfigure[Gradient norm $|\nabla \phi|(t)$ ($\Delta x = 1.25 \cdot 10^{-3}$).]{\includegraphics[width=0.5\columnwidth]{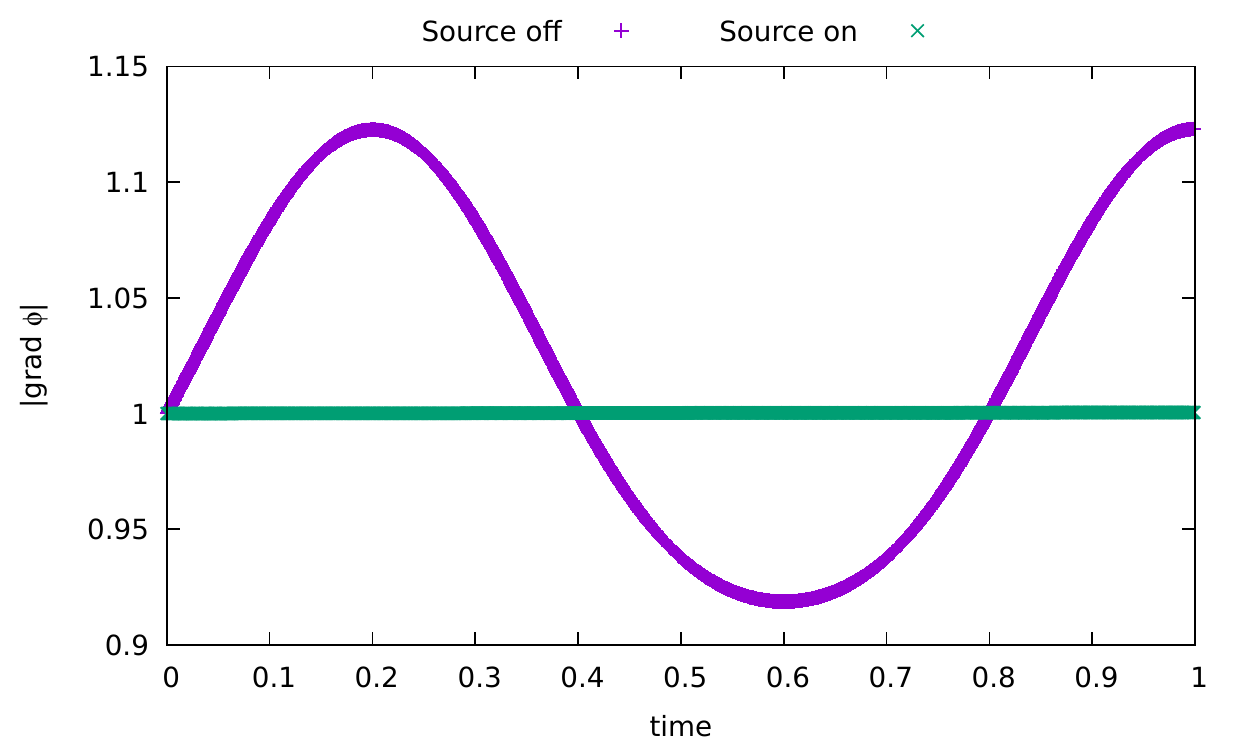}}
\subfigure[Mesh study for $|1-|\nabla \phi|(t)|$.]{\includegraphics[width=0.5\columnwidth]{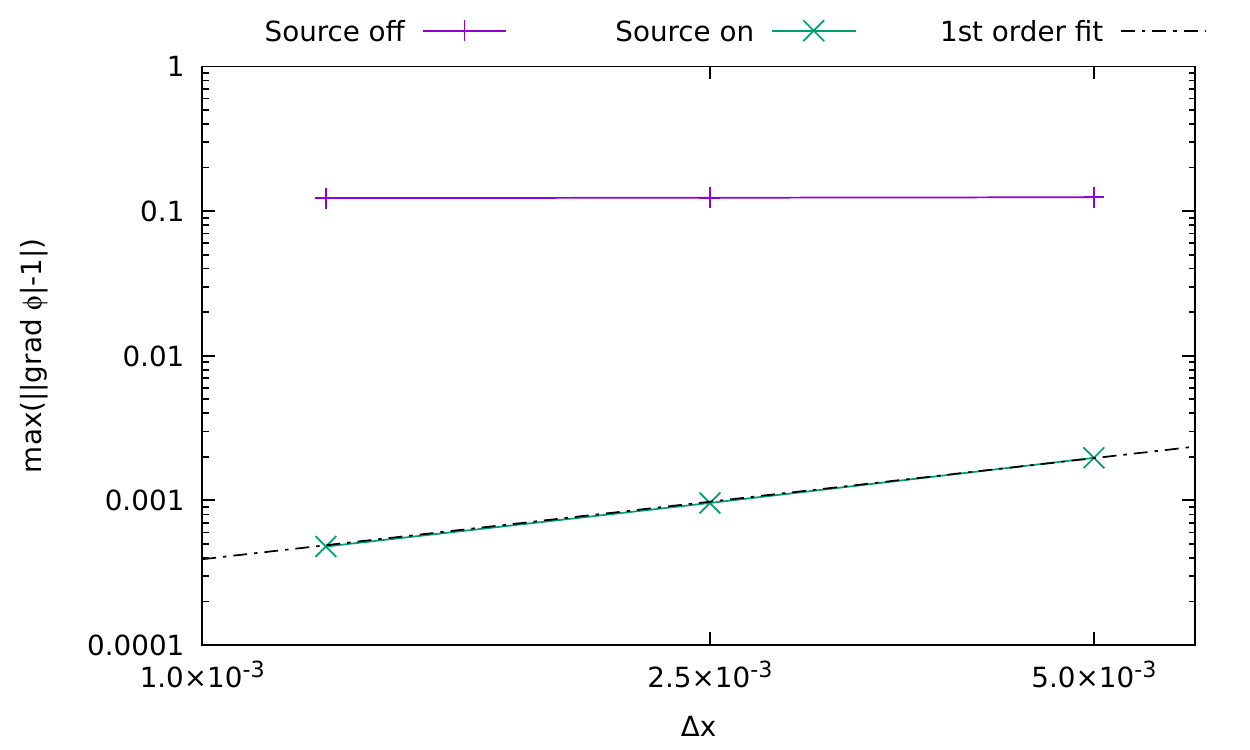}}
\caption{Numerical transport of the contact angle, curvature, and gradient norm for the time-periodic field \eqref{eqn:periodic_field}.}
\label{fig:periodic/position_and_angle_source_on_and_off}
\end{figure}

\paragraph{Numerical results:} The evolution of the interface for the ``Vortex-in-a-box field`` \eqref{eqn:shear_field} is shown in Figure~\ref{fig:dynamics-shear-2d}. The right contact point is moving further to the right while the contact angle approaches 90 degrees and the curvature at the contact point is decreasing. Figure~\ref{fig:shear/convergence_study} compares the numerical transport of the position, the contact angle, and the curvature at the contact point quantitatively. We observe that the results with and without source term are practically undistinguishable for the position, the contact angle and the curvature.
The mesh study for the maximum errors 
\begin{align*}
\max_{n} |x(t_n) - x_\text{ref}(t_n)|, \ \max_{n} |\theta(t_n) - \theta_\text{ref}(t_n)|, \ \max_{n} |\kappa(t_n) - \kappa_\text{ref}(t_n)|,
\end{align*}
reported in Figure~\ref{fig:shear/convergence_study} shows that all these quantities are converging with first-order. This is true regardless if the source term is active or not.
On the finest mesh, the maximum error for the contact angle is below $0.1$ degrees over the considered simulation time. Hence, indeed, the order of accuracy is preserved for the method with active source term.\\
\\
A major difference is, however, observed when we study the numerical evolution of $|\nabla \phi|$ as reported in Figures~\ref{fig:norm_grad_phi_field_plot} and \ref{fig:shear/gradient_norm_evolution}. \change{Figure~\ref{fig:norm_grad_phi_field_plot} shows a snapshot of $|\nabla \phi|$ close to the zero contour at time $t=0.875$. It is found that the gradient norm stays close to one in a neighborhood of the zero contour if the source term is active. This neighborhood corresponds to the region where the mollifier function is equal to one, i.e.\ to the region where $-w_1 \leq \phi \leq w_1.$ As expected, the conservation of the gradient norm works both away from the boundary and at the contact line. In contrast to that, a strong variation of $|\nabla \phi|$ along the zero contour is found for the standard level set method. In fact, in this case, the gradient norm is not conserved but evolving according to equation~\eqref{eqn:change-of-norm-standard-levelset}.}\\
\\
\change{Figure~\ref{fig:norm_grad_phi_field_plot} quantifies the evolution of $|\nabla \phi|$ along a flow trajectory following the contact line.} While $|\nabla \phi(t)|$ converges to a monotonically growing function (which is determined by equation~\eqref{eqn:change-of-norm-standard-levelset}) for the standard upwind scheme, we find a converge towards unity for the new method with active source term. By plotting the maximum deviation from a local-signed distance function, i.e.\
\[ \max_{n} |1-|\nabla \phi(t_n)||, \]
we see that the solution converges to a local signed-distance function with first-order for the adapted method with active source term. This is to be expected because of the time explicit approximation of the source term. The results reported in Figure~\ref{fig:periodic/position_and_angle_source_on_and_off} show that the same is true for the time-periodic example \eqref{eqn:periodic_field}.In summary, the simple adapted upwind method \eqref{eqn:adapted_first_order_upwind} works as expected.

\clearpage
\subsection{Numerical example in three dimensions}
\begin{figure}[hb]
\centering 
\includegraphics[width=0.6\columnwidth]{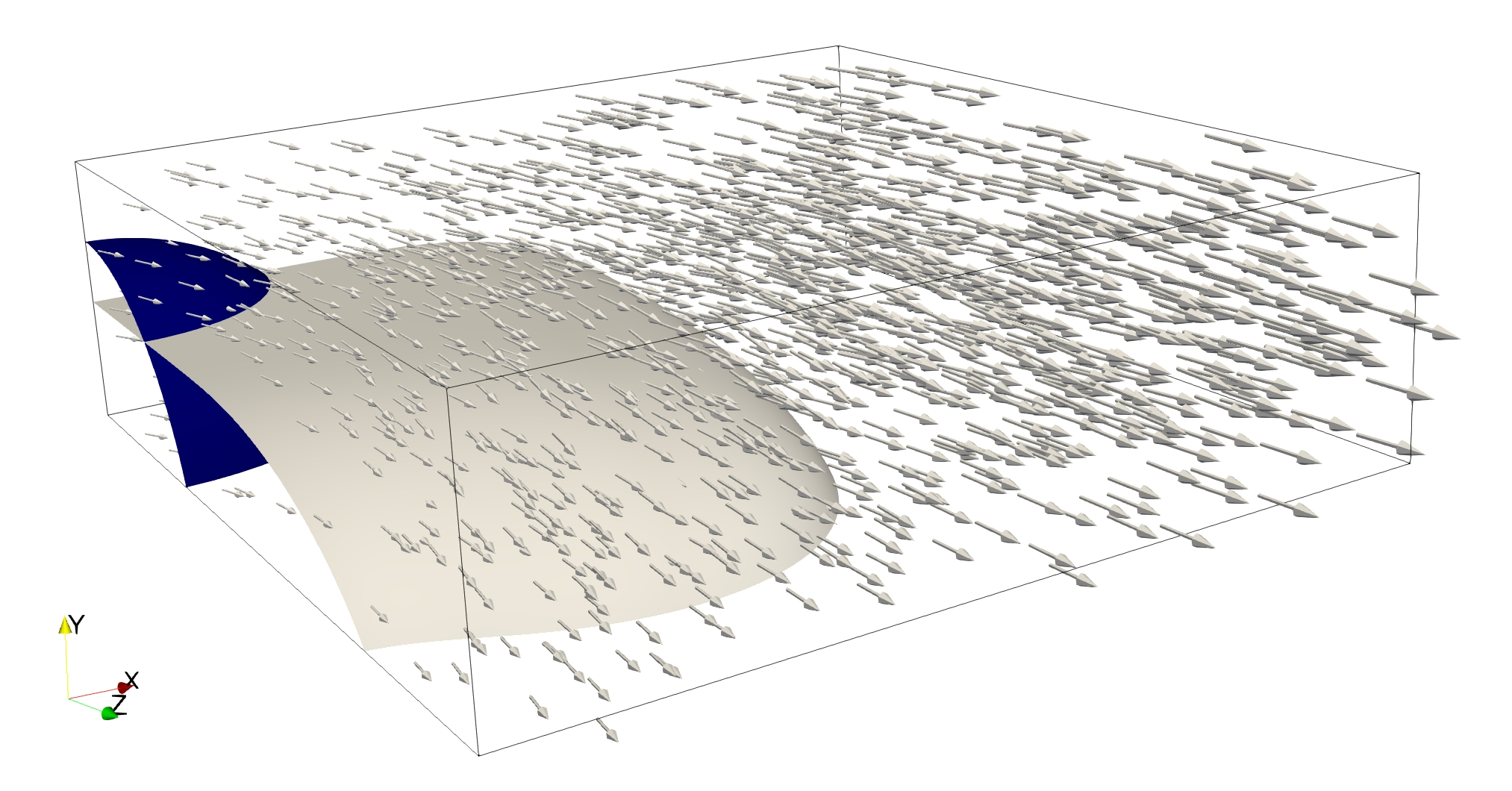}
\caption{Interface at $t=0$ and $t=1.93$ for the 3D advection test case.}\label{fig:3d_advection_test}
\end{figure}

We now consider a three-dimensional domain
\[ 0 \leq x \leq 2, \quad 0 \leq y \leq 0.6, \quad 0 \leq z \leq 2 \]
and choose the initial level set function
\[ \change{\phi_0(x,y) = \sqrt{(x-x_0)^2 + (y-y_0)^2 + (z-z_0)^2} - R_0 } \]
for $x_0 = z_0 = 0$, $y_0=-0.2$ and $R_0=0.6$. We follow the point on the contact line initially located at $(0.4,0,0.4)$. The initial interface normal at this point is $\nsigma = (2,1,2)/3$ corresponding to an initial contact angle $\theta_0 = \arccos(1/3)\approx 70.5^\circ$. The initial curvature of is $\kappa_0=-2/R = -10/3$. We choose a velocity field with Cartesian components $(v_1,v_2,v_3)$ given as
\begin{equation}
\begin{aligned}
v_1(x,y,z) &= v_1^0 + c_1 x + c_2 y + c_3 z,\\
v_2(x,y,z) &= -(c_1 + c_6) y,\\
v_3(x,y,z) &= v_3^0 + c_4 x + c_5 y + c_6 z.
\end{aligned}
\end{equation}
Note that $v$ satisfies incompressibility everywhere and impermeability at $y=0$ (i.e., $v_2(x,0,z)=0$) for any choice of parameters $v_1^0, v_3^0, c_1, c_2, c_3, c_4, c_5$ and $c_6 \in \RR$. We choose the following set of parameters
\[ v_1^0 = 0.3, \ v_3^0 = 0.4, \ c_1 = c_2 = - c_5 = c_6 = 0.1, \ c_3 = -0.2, \ c_4 = 0.3. \]
Figure~\ref{fig:3d_advection_test} shows the evolution of the interface and a visualization of the velocity field. Please note that we impose homogeneous Neumann boundary conditions for the level set function at inflow boundaries, i.e.\ at boundaries where the inequality $v \cdot \ndomega < 0$ holds. We choose an equidistant Cartesian mesh with $\Delta x = \Delta y = \Delta z$. For the mesh study presented below, the mesh resolution takes the values $\Delta x = 2/50, 2/100, 2/200$ corresponding to $15$ to $60$ cells per initial radius of curvature. The CFL number is kept constant at $0.2$.

\begin{figure}[h]
\subfigure[Results for the contact angle ($\Delta x = 0.01$).]{\includegraphics[width=0.5\columnwidth]{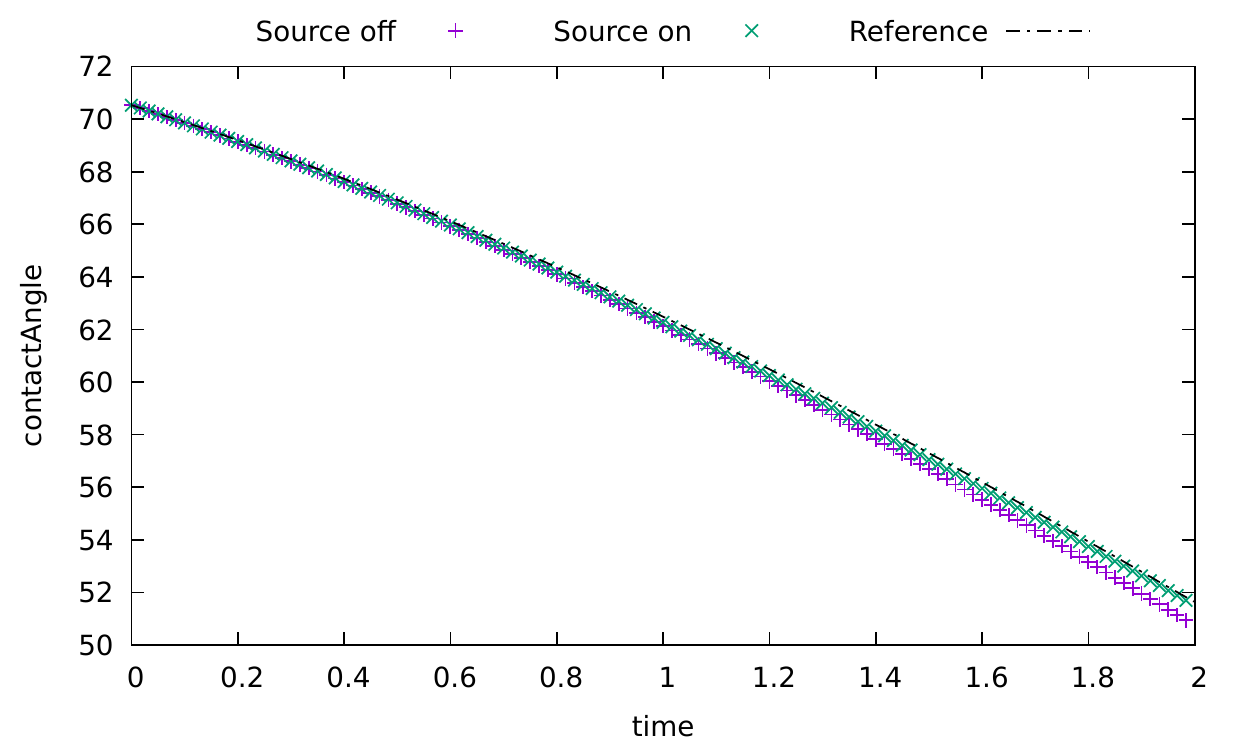}} 
\subfigure[Mesh study for the contact angle.]{\includegraphics[width=0.5\columnwidth]{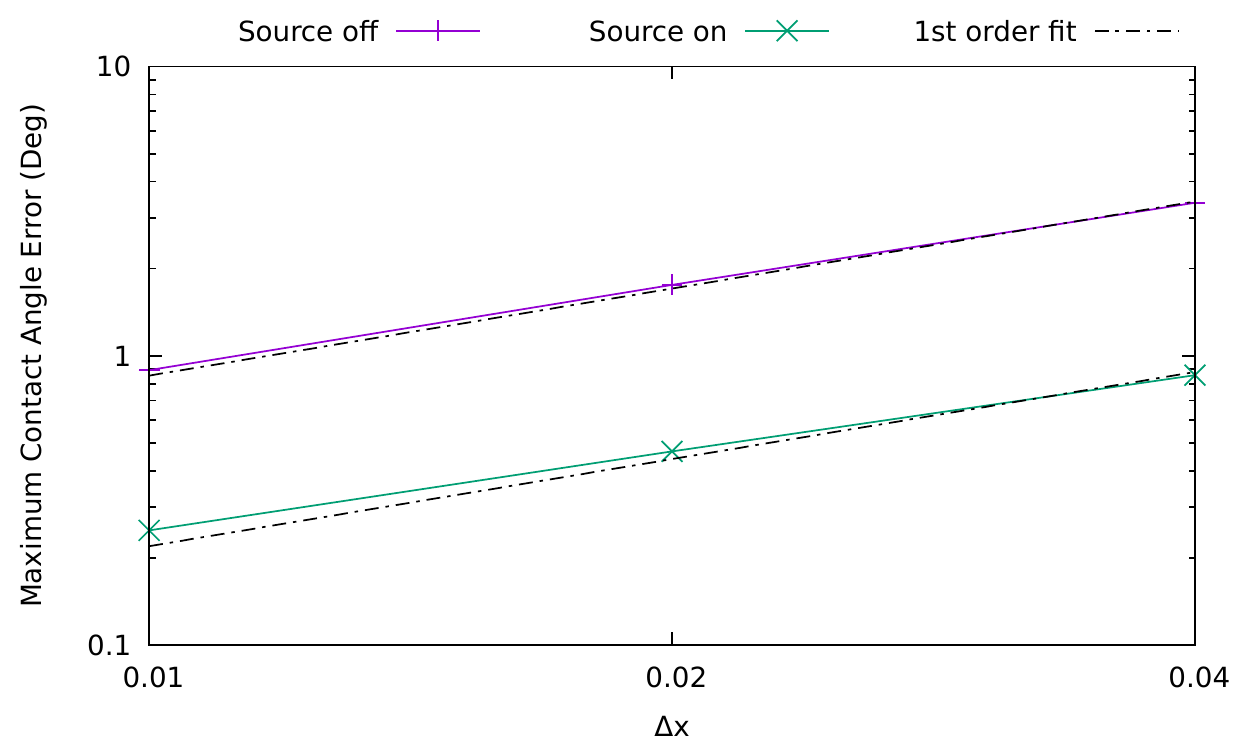}} 
\subfigure[Results for the mean curvature ($\Delta x = 0.01$).]{\includegraphics[width=0.5\columnwidth]{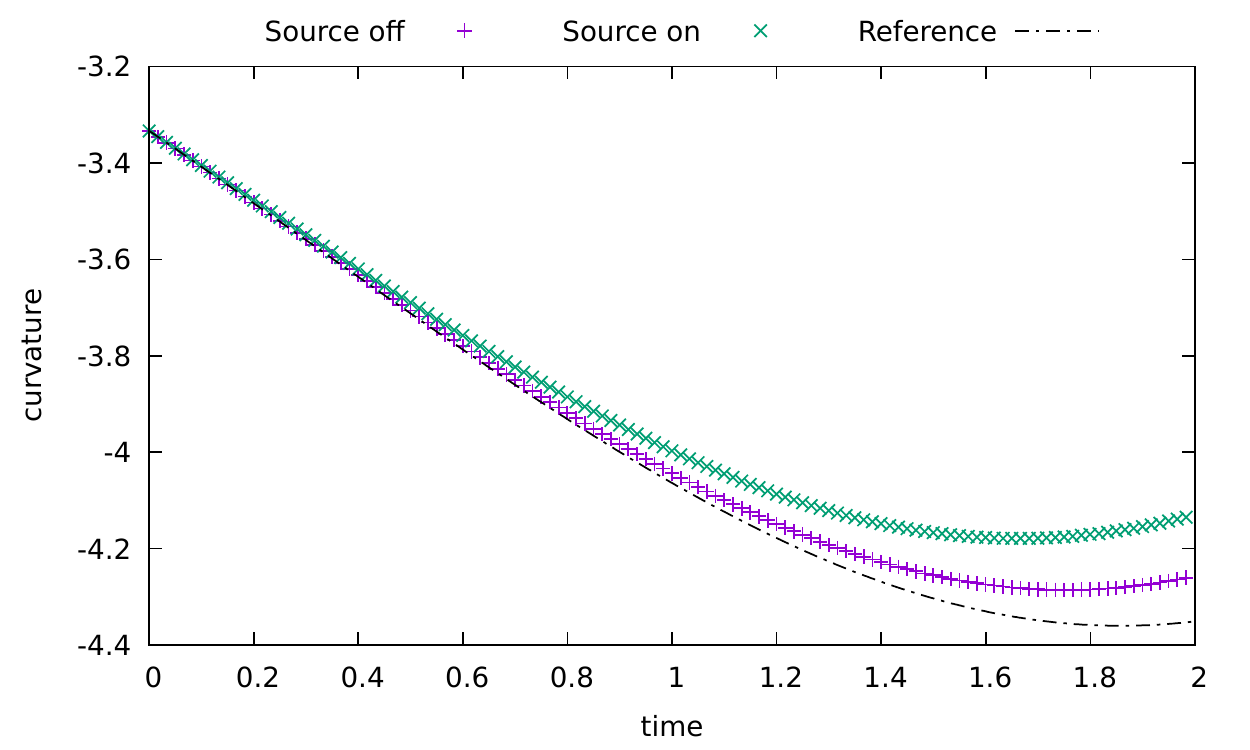}} 
\subfigure[Mesh study for the mean curvature.]{\includegraphics[width=0.5\columnwidth]{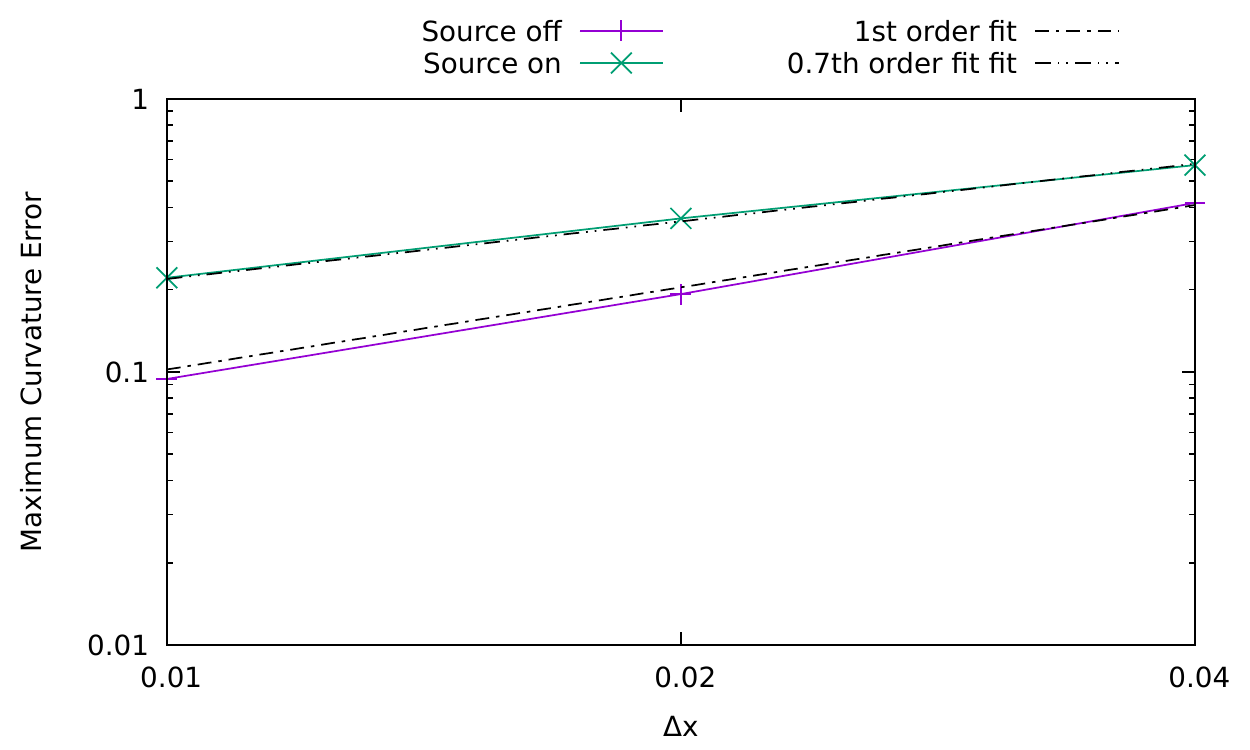}} 
\caption{Numerical evolution of the contact angle and the mean curvature in 3D.}\label{fig:3d_angle_and_curvature}
\end{figure}

\begin{figure}[h]
\subfigure[Gradient norm $|\nabla \phi|(t)$ ($\Delta x = 0.01$).]{\includegraphics[width=0.5\columnwidth]{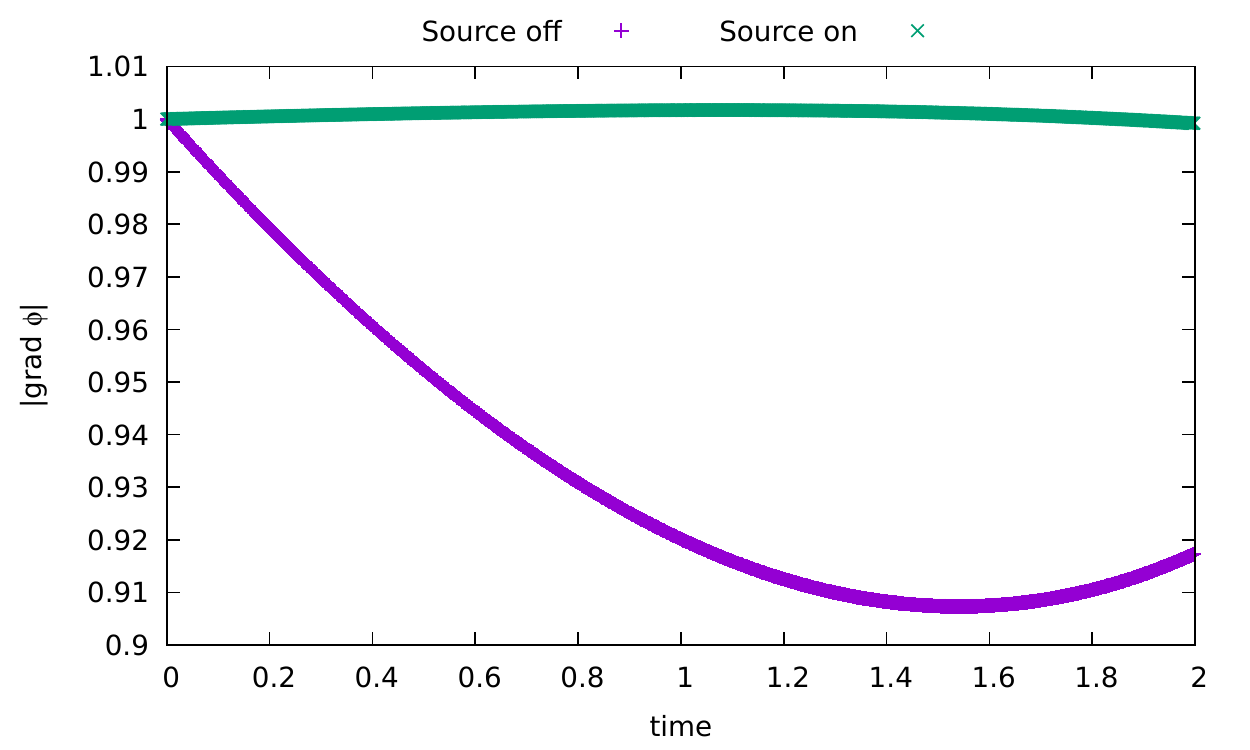}} 
\subfigure[Mesh study for $|1-|\nabla \phi|(t)|$.]{\includegraphics[width=0.5\columnwidth]{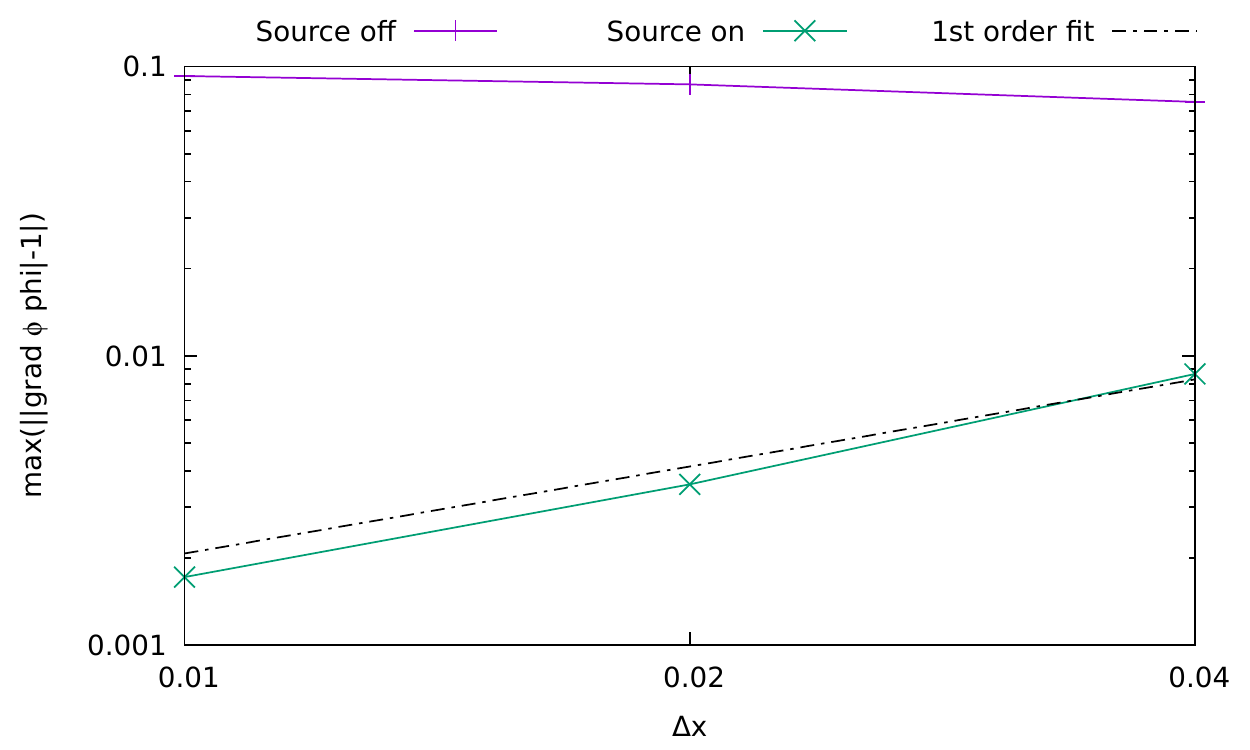}} 
\caption{Numerical evolution of $|\nabla \phi(t)|$ in 3D.}\label{fig:3d_gradient_norm}
\end{figure}

\paragraph{Results:} The results presented in Figure~\ref{fig:3d_angle_and_curvature} show that the transported contact angle converges to the reference solution with first-order and the accuracy for the contact angle is (in this case) improved with active source term. On the other hand, the accuracy for the mean curvature is decreased by source term. Notably, in this case, the order of convergence degenerates slightly from first-order with active source term. As expected, the gradient norm $|\nabla \phi(t)|$ converges to the constant 1 with first-order if the source term is active; see Figure~\ref{fig:3d_gradient_norm}. In summary, the adapted upwind method works as expected also in the three dimensional case (with some deviations for the curvature transport).

\section{Conclusion and outlook}
To summarize, we have introduced a modified level set transport equation that intrinsically preserves the norm of the gradient, i.e.\ $|\nabla \phi|$, at the physical interface represented by the zero contour of $\phi$. This is achieved by introducing a source term on the right-hand side which, however, turns the linear hyperbolic problem into a non-linear one. It is important to note that the source term is proportional to the \emph{physical} rate of interface generation.\\
\\
In practice, one can avoid solving the non-linear equation \eqref{eqn:adapted_levelset_equation} by approximating the interface generation by its snapshot at the previous time-step. This way, the problem is approximated by a linear hyperbolic problem of the form
\begin{align}\label{eqn:simplified_problem}
\partial_t \phi + v \cdot \nabla \phi = -r(x) \phi.
\end{align}
Equation \eqref{eqn:simplified_problem} has a number of convenient properties. In particular, the zero contours of \eqref{eqn:simplified_problem} and the original level set equation ($r=0$) are exactly the same. Hence, no source of error for the interface motion is introduced on the continuous level. Moreover, the problem \eqref{eqn:simplified_problem} is well-posed in the presence of a moving contact line without the need to impose a boundary condition for $\phi$. Hence, the signed distance preserving level set method works without special treatment in the case of a moving contact line. This is particularly interesting because some conventional redistancing methods are showing difficulties in the presence of a contact line. Finally, equation \eqref{eqn:simplified_problem} can readily be solved with standard numerical methods. We demonstrated this with the first-order upwind method which requires very little implementation work. First numerical examples show that the accuracy of the original upwind method is retained in terms of the interface position, the contact angle as well as the curvature, while the norm of the gradient at the interface converges to a constant value over time as the mesh is refined.\\
\\
Altogether, a first proof of concept for the locally signed distance preserving level set method (SDPLS) is established. Starting from here, we will work on an implementation in higher-order level set method and consider the coupling to a two-phase Navier Stokes solver. Among other things, it will be interesting to study if the method has any impact on the numerical conservation of the phase volume.
 
\paragraph{Acknowledgements:} Funded by the Deutsche Forschungsgemeinschaft (DFG, German Research Foundation) – Project-ID 265191195 – SFB 1194.\\
The authors gratefully acknowledge the computing time provided to them on the high-performance computer Lichtenberg at the NHR Centers NHR4CES at TU Darmstadt. This is funded by the Federal Ministry of Education and Research, and the state governments participating on the basis of the resolutions of the GWK for national high performance computing at universities (www.nhr-verein.de/unsere-partner).\newline
We thank Prof.\ Kohei Soga (Keio University, Tokio) for valuable discussion on the method.

\appendix
\section{Preliminaries on (moving) interfaces}\label{sec:appendix_interfaces_and_kinematics}
We briefly recall some basic mathematical definitions for (moving) hypersurfaces; see \cite{Giga.2006,Kimura2008,Pruess2016} for more details.

\paragraph{Construction of the signed distance function (see Chapter 2.3 in \cite{Pruess2016}):}
Mathematically, the signed distance function $d_\Sigma$ is constructed by inverting the map
\[ h \mapsto x + h \nsigma(x) \]
for some point $x$ on the interface $\Sigma$.

\begin{lemma}
Let $\Sigma \subset \RR^d$ be a $\mathcal{C}^2$-hypersurface and $x_0$ be an inner point of $\Sigma$. Then there exists an open neighborhood $U \subset \RR^n$ of $x_0$ and $\varepsilon > 0$ such that the map
$$
\begin{aligned}
&X:(\Sigma \cap U) \times(-\varepsilon, \varepsilon) \rightarrow \mathbb{R}^{n} \\
&X(x, h):=x+h n_{\Sigma}(x)
\end{aligned}
$$
is a diffeomorphism onto its image
$$
\mathcal{N}^{\varepsilon}:=X((\Sigma \cap U) \times(-\varepsilon \times \varepsilon)) \subset \mathbb{R}^{n} \text {, }
$$
i.e. $X$ is invertible there and both $X$ and $X^{-1}$ are $\mathcal{C}^{1}$. The inverse function has the form
\begin{align}\label{eqn:signed_distance_definition}
X^{-1}(x)=(\pi_\Sigma(x), d_\Sigma(x))
\end{align}
with $\mathcal{C}^{1}$-functions $\pi_\Sigma$ and $d_\Sigma$ on $\mathcal{N}^{\varepsilon}$.
\end{lemma}
In fact, the signed distance function $d_\Sigma$ is defined only locally on the so-called ``tubular neighborhood'' $\mathcal{N}^{\varepsilon}$. The operator $\pi_\Sigma$ in \eqref{eqn:signed_distance_definition} is the projection operator that maps each point $x \in \mathcal{N}^{\varepsilon}$ to its ``base point'' $\pi_\Sigma(x)$ located at the signed distance $d_\Sigma(x)$. Moreover, one can show that \cite{Pruess2016}
 \begin{align}
 \nabla d_\Sigma|_\Sigma = \nsigma
 \end{align}
 which implies that
 \begin{align}
 |\nabla d_\Sigma| = |\nsigma| = 1 \quad \text{on} \quad \Sigma.
 \end{align}
 As a consequence, the mean curvature of the interface can simply be computed as
 \begin{align}
 \kappa = - \Delta d_\Sigma \quad \text{at} \quad \Sigma.
 \end{align}
 
\paragraph{Moving hypersurfaces:} Similar definitions of a moving hypersurface can also be found in \cite{Giga.2006,Kimura2008,Pruess2016}. A generalization to moving hypersurfaces with boundary is given in \cite{Fricke2019}.
\begin{enumerate}[(a)]
 \item Let $I=(a,b)$ be an open interval. A family $\{\Sigma(t)\}_{t \in I}$ with $\Sigma(t) \subset \RR^3$ is called a \emph{$\mathcal{C}^{k,m}$-family of moving hypersurfaces} if the following holds.
\begin{enumerate}[(i)]
 \item Each $\Sigma(t)$ is an orientable $\mathcal{C}^m$-hypersurface in $\RR^3$ with unit normal field denoted as $\nsigma(t,\cdot)$.
 \item The graph of $\Sigma$, given as
 \begin{align}
 \label{eqn:def_moving_interface}
 \mathcal{M}:= \gr \Sigma = \bigcup_{t \in I} \{t\} \times \Sigma(t) \subset \RR\times\RR^{3},
 \end{align}
 is a $\mathcal{C}^k$-hypersurface in $\RR \times \RR^3$.
 \item The unit normal field is $k$-times continuously differentiable on $\mathcal{M}$, i.e.
 \[ \nsigma \in \mathcal{C}^k(\mathcal{M}). \]
\end{enumerate}
 \item Let $x_0 \in \Sigma(t_0)$ and $\gamma: I \rightarrow \mathbb{R}^{3}$ be a $\mathcal{C}^{1}$-curve on $\operatorname{gr} \Sigma$ that passes through $x_0$, i.e.
\[ \gamma(t_0) = x_0 \quad \text{and} \quad \left(t, \gamma(t)\right) \in \operatorname{gr} \Sigma \quad \forall t \in I. \]
Then, the \emph{normal speed} of $\Sigma(t_0)$ at $x_0$ is defined as
\begin{align} 
\normalspeed\left(t_{0}, x_{0}\right):= \gamma'(t_0) \cdot n_{\Sigma}\left(t_{0}, x_{0}\right) \in \RR. 
\end{align}
Note that the normal speed defined above does not depend on the choice of $\gamma$ (see \cite{Pruess2016} for details).
\end{enumerate}

\section{Surface transport theorem}\label{section:surface_transport}
We briefly recall the surface transport theorem for two-phase flow (see, e.g, \cite{Aris1989,Bothe2005}). Let $\{\Sigma(t)\}_{t \in I}$ be a family of moving surfaces and $\mathbf{v}^{\Sigma}: \operatorname{gr}(\Sigma) \rightarrow \mathbb{R}^{3}$ a consistent surface velocity field, i.e.\ $\mathbf{v}^{\Sigma}$ satisfies
\[ v^\Sigma \cdot \nsigma = \normalspeed \quad \text{on} \ \gr\Sigma. \]
Let $\{A(t)\}_{t \in I}$ be a co-moving control area inside $\Sigma(\cdot)$ and $\Phi^{\Sigma}: \operatorname{gr}(\Sigma) \rightarrow \mathbb{R}$ such that $\frac{\mathrm{D}^{\Sigma} \phi^{\Sigma}}{\mathrm{D} t}$ exists. Then
$$
\frac{\mathrm{d}}{\mathrm{d} t} \int_{A(t)} \phi^{\Sigma}(t, \mathbf{x}) \mathrm{d} o=\int_{A(t)}\left(\frac{\mathrm{D}^{\Sigma} \phi^{\Sigma}}{\mathrm{D} t}(t, \mathbf{x})+\phi^{\Sigma} \operatorname{div}_{\Sigma} \mathbf{v}^{\Sigma}(t, \mathbf{x})\right) \mathrm{d} o
$$
This yields the following relation for co-moving areas
$$
\frac{\mathrm{d}}{\mathrm{d} t}|A(t)|=\frac{\mathrm{d}}{\mathrm{d} t} \int_{A(t)} \mathrm{d} o=\int_{A(t)} \operatorname{div}_{\Sigma} \mathbf{v}^{\Sigma} \mathrm{d} o.
$$

\end{document}